\documentclass[12pt, twoside]{article}
\usepackage[english]{babel}
\usepackage{amssymb}
\usepackage{amsmath}
\usepackage{epsfig}
\let\noi\noindent

\newcommand{\an}{\mbox{\scriptsize an}}
\newcommand{\fin}{\mbox{\scriptsize fin}}

\def\exp{\text{exp}}

\def\Arg{\text{Arg}}
\def\Sing{\text{Sing}}

\def\supp{\text{supp}}
\def\sup{\text{sup}}

\def\Re{\text{Re}}

\def\id{\text{id}}
\def\Im{\text{Im}}

\begin{document}

$\phantom{^{*)}}$ \insert\footins{\vskip3pt
\setbox0\hbox{$\phantom{^{*)}}$}\hangindent\wd0 $\phantom{^{*)}}$
{{\footnotesize
\hspace{-0.78cm} {\it 2000 Mathematics Subject Classification}. 03C64, 32B20, 35C20, 35J25,
30E15, 30D60, 30D05, 37E35\\
{\it Keywords and phrases}: o-minimal structures, semianalytic sets, Dirichlet problem,
quasianalytic classes, dynamical systems\\
Supported by DFG-project KN202/5-2 }}}

\noindent {\large  THE DIRICHLET PROBLEM IN THE PLANE WITH \break
SEMIANALYTIC RAW DATA, QUASIANALYTICITY AND\break O-MINIMAL
STRUCTURES}

\noindent $\hbox to 14cm{\hrulefill}$

\medskip\noi
TOBIAS KAISER

\vspace{0.6cm}\noi {\bf Abstract}

\noindent {\it We investigate the Dirichlet solution for a
semianalytic continuous function on the boundary of a semianalytic
bounded domain in the plane. We show that the germ of the Dirichlet
solution at a boundary point with angle greater than 0 lies in a
certain quasianalytic class used by Ilyashenko in his work on
Hilbert's 16$^{\mbox{\tiny th}}$ problem. With this result we can
prove that the Dirichlet solution is definable in an o-minimal
structure if the angle at a singular boundary point of the domain is
an irrational multiple of $\pi$.}

\vspace{1.2cm}\noi {\bf Introduction}

\medskip\noi
Traditional and excellent settings for `tame geometry' on the reals
are given by the category of semialgebraic sets and functions and by
the category of subanalytic sets and functions. The sets considered
may have singularities but behave still `tame', i.e. various
finiteness properties hold, see Bierstone-Milman [4], Bochnak et al.
[5], Denef-Van den Dries [8], \L ojasiewicz [32] and Shiota [35].
These categories are excellent for geometrical questions but as
often observed they are insufficient for problems from analysis. For
example, the solution of the differential equation $y' =
\frac{y}{x^2}$ on ${\mathbb R}_{> 0}$, given by $x \mapsto
e^{-\frac{1}{x}}$, is not subanalytic anymore. Therefore a natural
aim was a better understanding of the solutions of first order
ordinary differential equations or more general of Pfaffian
equations with polynomial or analytic raw data, and a lot of
research activities was done in this direction.

\noi It was shown that sets defined by the solutions of Pfaffian
equations, so-called semi- and sub-Pfaffian sets, show a `tame'
behaviour, see for example Cano et al. [6], Gabrielov [18],
Gabrielov et al. [19] and Lion-Rolin [30]. Also a more axiomatic
understanding was obtained. This axiomatic setting is given by the
framework of o-minimal structures. They generalize the category of
semialgebraic sets and functions and are defined by finiteness
properties. They are considered as ``an excellent framework for
developing tame topology, or  {\it topologie mod\'er\'ee}, as
outlined in Grothendieck's prophetic ``Esquisse d'un Programme'' of
1984'' (see the preface of Van den Dries [10], which provides a very
good source for the definition and the basic properties of o-minimal
structures). The basic example for an o-minimal structure is given
by the semialgebraic sets and functions; these are the sets which
are definable from the real field ${\mathbb R}$ by addition,
multiplication and the order. The subanalytic category fits not
exactly in this concept (compare with Van den Dries [9]), but the
globally subanalytic sets, i.e. the sets which are subanalytic in
the ambient projective space, form an o-minimal structure, denoted
by ${\mathbb R}_{\an}$ (see Van den Dries-Miller [13]). A
breakthrough was achieved by Wilkie, who showed in [40], using
Khovanskii theory for Pfaffian systems (see [29]), that the real
exponential field ${\mathbb R}_{\exp}$, i.e. the field of reals
augmented with the global exponential function $\exp \colon {\mathbb
R} \to {\mathbb R}_{> 0}$, is an o-minimal structure. Subsequently
Van den Dries-Miller [12] and Van den Dries et al. [11] proved that
the structure ${\mathbb R}_{\an, \exp}$ is o-minimal. For general
Pfaffian functions o-minimality was again obtained by Wilkie [41].
This result was extended by Karpinski-Macintyre [28] and finally
stated by Speissegger [36] in its most generality: the Pfaffian
closure of an o-minimal structure on the real field is again
o-minimal.

\noi So first order differential equations or more general Pfaffian
equations in the subanalytic context resp. in the context of
o-minimal structures are well understood. As an application
integration of a one variable function in an o-minimal structure can
be handled (for integration with parameters this is the case so far
only for subanalytic maps by the results of Lion-Rolin [31] and
Comte~et~al.~[7]).

\noi Our goal is to attack higher order partial differential
equations in the subanalytic resp. o-minimal setting. Compared to
ordinary differential equations there are distinct classes of
equations and boundary value problems, each with its own theory. A
very important class of PDE's is given by the elliptic ones and one
of its outstanding representative is the Laplace equation. We
consider the Dirichlet problem, i.e. the Laplace equation with
boundary value problem of the first kind: let $\Omega \subset
{\mathbb R}^n$ be a bounded domain and let $h \in C(\partial
\Omega)$ be a continuous function on the boundary. Then the
Dirichlet problem for $h$ is the following: is there a function $u$
continuous on $\overline{\Omega}$ and twice differentiable in
$\Omega$ such that
$$ \begin{array}{rcl}
\Delta u & = & 0 \quad \mbox{in} \quad \Omega,\\
u & = & h \quad \mbox{on} \quad \partial \Omega.\\
\end{array}$$

\noi Here $\Delta := \frac{\partial^2}{\partial x_1^2} + \dots +
\frac{\partial^2}{\partial x_n^2}$ is the Laplace operator.
Functions fulfiling the first equality are called harmonic in
$\Omega$. They are actually real analytic on $\Omega$. If the answer
is yes, i.e. if such a $u$ exists, we call it the Dirichlet solution
for $h$. If the answer is yes for all continuous boundary functions
the domain $\Omega$ is called regular. The punctured open ball for
example is irregular (see Helms [20, p.168]). Simply connected
domains in the plane are regular. In [25] we gave in the case that
$\Omega$ is subanalytic a necessary and sufficient condition for
$\Omega$ to be regular. For irregular domains there is a more
general solution for the Dirichlet problem, the so-called
Perron-Wiener-Brelot solution (see Armitage-Gardiner [2, Chapter 6]
and [20, Chapter 8]).\\
We are interested in the case that $\Omega$ is a subanalytic domain
and that also the boundary function $h$ is subanalytic. The natural
questions are now the following. What can be said about the
Dirichlet solution? Is it definable in an o-minimal structure? We
consider the case that $\Omega$ is a domain in the plane (then
$\Omega$ and $h$ are semianalytic; see [4, Theorem 6.13]). By [25]
$\Omega$ is regular if it has no isolated boundary points. If this
is not the case the Perron-Wiener-Brelot solution for a continuous
boundary function coincides with the Dirichlet solution for this
function after adding the finitely many isolated boundary points to
$\Omega$ (compare with [20, p.168]). So from now on we may assume
that $\Omega$ has no isolated boundary points. Under the additional
assumption that the boundary is analytically smooth it was shown in
[24] that the Dirichlet solution is definable in the o-minimal
structure ${\mathbb R}_{\an, \exp}$. This result is obtained by
reducing the problem to the unit ball. There the Dirichlet solutions
are given by the Poisson integral (see [2, Chapter 1.3]) and we can
apply the results about integration of subanalytic functions (see
[7] and [31]).

\noi The challenging part are domains with singularities. The
starting point to attack singularities are asymptotic expansions.
Given a simply connected domain $D$ in ${\mathbb R}^2$ which has an
analytic corner at $0 \in \partial D$ (i.e. the boundary at 0 is
given by two regular analytic curves which intersect in an angle
$\sphericalangle D$ greater than 0) and a continuous boundary
function $h$ which is given by power series on these analytic
curves, Wasow showed in [39] that the Dirichlet solution for $h$ is
the real part of a holomorphic function $f$ on $D$ which has an
asymptotic development at 0 of the following kind:
$$f(z) \sim \sum\limits_{n = 0}^{\infty} a_n P_n (\log z)
z^{\alpha_n} \quad \mbox{as} \quad z \longrightarrow 0 \quad
\mbox{on} \quad D, \leqno(\dagger)$$

\noi i.e. for each $N \in {\mathbb N}_0$ we have
$$f(z) - \sum\limits_{n = 0}^N a_n P_n (\log z)
z^{\alpha_n} = o (z^{\alpha_n}) \quad \mbox{as} \quad z
\longrightarrow 0 \quad \mbox{on} \quad D, $$

\noi where $\alpha_n \in {\mathbb R}_{\geq 0}$ with $\alpha_n
\nearrow \infty$, $P_n \in {\mathbb C} [z]$ monic and $a_n \in
{\mathbb C}^*$. Moreover, if $\sphericalangle D/\pi \in {\mathbb R}
\setminus {\mathbb Q}$ we have that $P_n = 1$ for all $n \in
{\mathbb N}_0$. (Note that $P_0 = 1$ for any angle).\\
To use this asymptotic development we want to have a quasianalytic
property; we want to realize these maps in a class of functions with
an asymptotic development as in ($\dagger$) such that the functions
in this class are determined by the asymptotic expansion. Such
quasianalyticity properties are key tools in generating o-minimal
structures (see [27], Van den Dries-Speissegger [14, 15] and Rolin
et al. [34]; see also Badalayan [3] for quasianalytic classes of
this kind).

\noi Exactly the same kind of asymptotic development occurs at the
transition map of a real analytic vector field on ${\mathbb R}^2$ at
a hyperbolic singularity (see Ilyashenko [22]). Poincar\'e return
maps are compositions of finitely many transition maps and are an
important tool to the qualitative understanding of the trajectories
and orbits of a polynomial or analytic vector field on the plane.
Following Dulac's approach (see [16]), Ilyashenko uses asymptotic
properties of the Poincar\'e maps to prove Dulac's problem (the weak
form of (the second part) of Hilbert's 16th problem): a polynomial
vector field on the plane has finitely many limit cycles (see
Ilyashenko [23] for an overview of the history of Hilbert~16, part
2). One of the first steps in Ilyashenko's proof is to show that the
transition map at a hyperbolic singularity is in a certain
quasianalytic class. Formulating his result on the Riemann surface
of the logarithm (compare with the introduction of [27] and with
[27, Proposition~2.8]) he proves that the considered transition maps
have a holomorphic extension to certain subsets of the Riemann
surface of the logarithm, so-called standard quadratic domains (see
[22, \S 0.3] and [27]; see also Section~1 below), such that the
asymptotic development holds there. Quasianalyticity follows then by
a Phragm\'en-Lindel\"of argument (see [22, \S 3.1]). The extension
of the transition map at a hyperbolic singularity and its asymptotic
development is obtained by transforming the vector field by a real
analytic change of coordinates into a suitable form. Then the
complexification of the resulting vector field is considered and
certain curves are lifted to the Riemann surfaces of its complex
phase curves (see [22, \S 0.3] and Ilyashenko [21, Section 3]). The
extension process can be rediscovered as a discrete dynamical system
given by the iteration of a local biholomorphic map which fixes the
origin and whose first derivative at the origin has absolute value 1
(see [21, Proposition 3]).

\medskip\noindent
We can link the Dirichlet problem with semianalytic raw data and
Ilyashenko's quasianalytic class:

\bigskip\noi
THEOREM A

\medskip\noi
{\it Let} $\Omega \subset {\mathbb R}^2$  {\it with} $0 \in
\partial \Omega$ {\it be a bounded
semianalytic domain without isolated boundary points. Let $h$ be a
semianalytic and continuous function on the boundary and let $u$ be
the Dirichlet solution for $h$. If the angle of $\Omega$ at $0$ is
greater than $0$, then $u$ is the real part of a holomorphic
function which is in the quasianalytic class of Ilyashenko described
above}.

\newpage\noi We will give a precise definition of 'the angle' of $\Omega$
at a boundary point in Section 2 below.

\medskip\noi
The semianalytic boundary curves of the given domain and the
semianalytic boundary function are locally given by Puiseux series
(see [9, p.192]). We obtain Theorem~A by extending the result of
Wasow to Puiseux series and by a geometric argument: we repeat
reflections of the Dirichlet solution at two real analytic curves to
go all the way up and down the Riemann surface of the logarithm to
get the extension to a standard quadratic domain. It is crucial that
the boundary function is given by Puiseux series to be able to apply
the reflection process. In the literature pairs of germs of real
analytic curves resp. groups generated by two non-commuting
antiholomorphic involutions (reflections at real analytic curves
correspond to antiholomorphic involutions) are studied in the
context of the classification of germs of biholomorphic maps fixing
the origin which is investigated by the Ecalle-Voronin theory (see
Ecalle [17], Voronin [37, 38] and Ahern-Gong [1]). Local
biholomorphic maps which fix the origin and whose first derivatives
at the origin have absolute value 1 occur also at the reflection
process: the description of the repeated reflections of the
Dirichlet solution involves iteration, inversion and conjugation of
such maps (resp. their lifting to the Riemann surface of the
logarithm). Also summation of Puiseux series is involved. We
carefully estimate the functions obtained by the reflection process
to get Theorem~A.

\bigskip\noi Transition maps at a hyperbolic singularity exhibit a
similar dichotomy of the asymp\-to\-tic development as indicated in
($\dagger$), depending on whether the hyperbolic singularity is
resonent or non-resonant, i.e. whether the ratio of the two
eigenvalues of the linear part of the vector field at the given
hyperbolic singularity is rational or irrational, see [16] and [27].
In [27] it is shown that transition maps at non-resonant hyperbolic
singularities are definable in a common o-minimal structure, denoted
by ${\mathbb R}_{{\mathcal Q}}$. This structure ${\mathbb
R}_{{\mathcal Q}}$ is generated by the functions (restricted to the
positive line) in Ilyashenko's quasianalytic class that have no
log-terms in their asymptotic expansion.

\medskip\noi
With this result and Theorem A, taking into account the nonsingular
boundary points (compare with [24]), we can prove

\newpage\noi THEOREM B

\medskip\noi
{\it Let} $\Omega \subset {\mathbb R}^2$ {\it be a bounded
semianalytic domain without isolated boundary points and let $h$ be
a semianalytic and continuous function on the boundary. Suppose that
the following condition holds: if $x$ is a singular boundary point
of $\Omega$ then the angle of the boundary at $x$ divided by $\pi$
is irrational. Then the Dirichlet solution for $u$ (i.e. its graph
considered as a subset of ${\mathbb R}^3$) is definable in the
o-minimal structure} ${\mathbb R}_{{\mathcal Q}, \exp}$.

\medskip\noi
As an application we obtain that the Green function of a bounded
semianalytic domain fulfiling the assumptions of Theorem B, is
definable in the o-minimal structure ${\mathbb R}_{{\mathcal Q},
\exp}$. If the considered domain is semilinear the assumption on the
angles can be dropped.

\medskip\noi
In [26] it is shown that the Riemann map from a simply connected
bounded and semianalytic domain in the plane with the same
assumptions on the angles as above to the unit ball is definable in
${\mathbb R}_{{\mathcal Q}}$. There it is also the key step to
realize the function in question in the quasianalytic class of
Ilyashenko. But the main ingredient, the reflection procedure,
differs heavily. There is some overlap in the definitions, but in
order to keep this paper reasonably self-contained we include all
necessary definitions here.

\medskip\noi
The paper is organized as follows: Section~1 is about the Riemann
surface of the logarithm and the classes of functions which we use
later. In particular Ilyashenko's quasianalytic class is introduced.
In Section~2 we define the notion of an angle for semianalytic
domains in a rigorous way and we present the concept of a domain
with analytic corner. In Section~3 we prove Theorem~A and Theorem~B
and give applications.

\bigskip\noi {\it Notation}

\medskip\noi
By ${\mathbb N}$ we denote the set of natural numbers and by
${\mathbb N}_0$ the set of nonnegative integers. Let $a \in {\mathbb
R}^n$ and $r > 0$. We set $B(a, r) := \{ z \in {\mathbb R}^n \;
\vert \; \vert z - a \vert < r \}$, $\overline{B} (a, r) := \{ z \in
{\mathbb R}^n \; \vert \; \vert z - a \vert \leq r \}$ and $\dot{B}
(a, r) := B(a, r) \setminus \{ a \}$. Here $\vert \quad \vert$
denotes the euclidean norm in ${\mathbb R}^n$. A domain is an open,
nonempty and connected set (in a topological space). Given an open
set $U$ of a Riemann surface we denote by ${\mathcal O} (U)$ the
${\mathbb C}$-algebra of holomorphic functions on $U$ with values in
${\mathbb C}$. By ${\mathcal O}_0$ we denote the ${\mathbb
C}$-algebra of germs of holomorphic functions in open neighbourhoods
of $0 \in {\mathbb C}$. Note that ${\mathcal O}_0 = {\mathbb C} \{ z
\}$. We identify ${\mathbb C}$ with ${\mathbb R}^2$.

\vspace{1.5cm}\noi {\bf 1. The Riemann surface of the logarithm}

\medskip\noi
We establish the setting on the Riemann surface of the logarithm. We
introduce the quasianalytic class of Ilyashenko, generalizing the
setting from [27, Section 2 \& 5] (we apply the results and methods
of [27] to obtain Theorem~B in Section 3). We define several other
classes of functions we will use in the proof of Theorem~A.

\medskip\noi
{\it Definition 1.1} (compare with [27, pp.12-13])

\noi We define the Riemann surface of the logarithm ${\mathbf L}$ in
polar coordinates by ${\mathbf L} := {\mathbb R}_{> 0} \times
{\mathbb R}$. Then ${\mathbf L}$ is a Riemann surface with the
isomorphic holomorphic projection map $\log \colon {\mathbf L} \to
{\mathbb C}$, $(r, \varphi) \mapsto \log r + i \varphi$. For $z =
(r, \varphi) \in {\mathbf L}$ we define the absolute value $\vert z
\vert := r$ and the argument arg $z := \varphi$. For $r > 0$ we set
$B_{{\mathbf L}}(r) := \{ z \in {\mathcal L}\; \vert \; \vert z
\vert < r \}$. We identify ${\mathbb C} \setminus {\mathbb R}_{\leq
0}$ with ${\mathbb R}_{> 0} \times ] - \pi, \pi [ \subset {\mathbf
L}$ via polar coordinates. Let $\alpha \geq 0$. We define the power
function $z^{\alpha}$ as $z^{\alpha} \colon {\mathbf L} \to {\mathbb
C}$, $z = (r, \varphi) \mapsto \exp (\alpha \log z)$.

\bigskip\noi {\it Remark 1.2}

\noi The logarithm on ${\mathbf L}$ extends the principal branch of
the logarithm on ${\mathbb C} \setminus {\mathbb R}_{\leq 0}$. The
power functions on ${\mathbf L}$ extend the power functions on
${\mathbb C} \setminus {\mathbb R}_{\leq 0}$.

\bigskip\noi
{\it Definition 1.3} (compare with [27, p.25])

\noi a) Let $\rho \geq 0$. We define the holomorphic map ${\mathbf
p}^{\rho} \colon {\mathbf L} \to {\mathbf L}$, $(r, \varphi)
\mapsto (r^{\rho}, \rho {\varphi})$.\\
b) We define the holomorphic map ${\mathbf m} \colon {\mathbf L}^2
\to {\mathbf L}$, $((r_1, \varphi_1), (r_2, \varphi_2)) \mapsto (r_1
r_2, \varphi_1 + \varphi_2)$.

\bigskip\noi {\it Remark 1.4}

\noi
\begin{itemize}
\item[a)] Let $r > 0$ and let $f \colon B (0, r) \to {\mathbb C}$ be
holomorphic. Then the function $f_{{\mathbf L}} \colon B_{{\mathbf
L}} (r) \to {\mathbb C}$, $z \mapsto f (z^1)$, is holomorphic
(compare with Definition 1.1).
\item[b)] Let $r > 0$ and let $g = \sum\limits_{n = 0}^{\infty} a_n
t^{\frac{n}{d}}$ be a convergent Puiseux series on $]0, r [$ (with
complex coefficients). Then $g_{{\mathbf L}} \colon B_{{\mathbf L}}
(r) \to {\mathbb C}$, $z \mapsto f_{{\mathbf L}} \left( {\mathbf
p}^{\frac{1}{d}} (z) \right)$, where $f(z) = \sum\limits_{n =
0}^{\infty} a_n z^n \in {\mathcal O} (B(0, r^d))$, is holomorphic.
Let $z \in B_{{\mathbf L}}(r)$. Then $g_{{\mathbf L}} (z) =
\sum\limits_{n = 0}^{\infty} a_n z^{\frac{n}{d}}$ as an absolutely
convergent sum, i.e. $\sum\limits_{n = 0}^{\infty} \vert a_n \vert
\; \vert z^{\frac{n}{d}} \vert < \infty$ (where $z^{\frac{n}{d}}$ as
defined in Definition 1.2). So we do often not distinguish between
$g$ and $g_{{\mathbf L}}$. Viewed as a function on ${\mathbf L}$ we
say that $g(z) = \sum\limits_{n = 0}^{\infty} a_n z^{\frac{n}{d}}$
is a Puiseux series convergent on $B_{{\mathbf L}}(r) \subset
{\mathbf L}$. If $d = 1$ (i.e. we are in case a)) we say that $g$ is
a power series convergent on $B_{{\mathbf L}}(r)$.
\end{itemize}

\bigskip\noi PROPOSITION 1.5

\noi {\it Let $r > 0$ and let $g(z) = \sum\limits_{n = 0}^{\infty}
a_n z^{\frac{n}{d}}$ be a Puiseux series convergent on $B_{{\mathbf
L}}(r) \subset {\mathbf L}$.}
\begin{itemize}
\item[a)] Let $t > r$ and assume that there is some $G \in
{\mathcal O} (B_{{\mathbf L}}(t))$ such that $G = g$ on $B_{{\mathbf
L}} (r)$. Then $\sum\limits_{n = 0}^{\infty} a_n z^{\frac{n}{d}}$ is
convergent on $B_{{\mathbf L}}(t)$.
\item[b)] {\it Let $c > 0$ such that $\vert g(z) \vert
\leq c$ for $\vert z \vert < r$. Then for every $N \in {\mathbb N}$
and every $\vert z \vert < r$, we obtain}

$$ \big\vert g (z) - \sum\limits_{n = 0}^N a_n
z^{\frac{n}{d}} \big\vert \leq c \left( \frac{\vert z \vert}{r}
\right)^{\frac{N+1}{d}} \left( \frac{1}{1- \left( \frac{\vert z
\vert}{r} \right)^{\frac{1}{d}}} \right).$$
\end{itemize}

\bigskip\noi
{\it Proof}

\noi
\begin{itemize}
\item[a)] We set $h \colon B_{{\mathbf L}} (t^{\frac{1}{d}}) \to
{\mathbb C}$, $z \mapsto G ({\mathbf p}^d (z))$. Then $h \in
{\mathcal O} (B_{{\mathbf L}}(t^{\frac{1}{d}}))$ and $h (z) =
\sum\limits_{n = 0}^{\infty} a_n z^n$ for $\vert z \vert <
r^{\frac{1}{d}}$. Let $a := (1, 2 \pi)$. Note that $h({\mathbf m}
(a, z)) = h (z)$ for all $\vert z \vert < r^{\frac{1}{d}}$ and hence
for all $\vert z \vert < t^{\frac{1}{d}}$. Hence $\tilde{h} \colon
\dot{B} (0, t^{\frac{1}{d}}) \to {\mathbb C}$, $z \mapsto h (\vert z
\vert e^{i \Arg z})$ (with $\Arg z \in ] - \pi, \pi]$ the standard
argument of a complex number $z \in {\mathbb C}^*$), is well defined
and holomorphic. We have $\tilde{h} (z) = \sum\limits_{n =
0}^{\infty} a_n z^n$ for $\vert z \vert < r^{\frac{1}{d}}$ and hence
for $\vert z \vert < t^{\frac{1}{d}}$ by Cauchy's Theorem. By Remark
1.4 we see that $G(z) = \tilde{h}_{{\mathbf L}} ({\mathbf
p}^{\frac{1}{d}} (z))$ is a Puiseux series which is convergent on
$B_{{\mathbf L}}(t)$.
\item[b)] Let $h (z) = \sum\limits_{n = 0}^{\infty} a_n
z^{\frac{n}{d}}$ be the corresponding Puiseux series on $B(0, r)$,
i.e. $g = h_{{\mathbf L}}$. Let $f(z) = \sum\limits_{n = 0}^{\infty}
a_n z^n \in {\mathcal O} (B(0, r^{\frac{1}{d}}))$. From the Cauchy
estimates we obtain $\big\vert f(z) - \sum\limits_{n = 0}^N a_n z^N
\big\vert \leq c \left( \frac{\vert z \vert}{r^{\frac{1}{d}}}
\right)^{N+1} \left( \frac{1}{1- \frac{\vert z
\vert}{r^{\frac{1}{d}}}} \right)$ and the claim follows, since $g =
f_{{\mathbf L}} \circ {\mathbf p}^{\frac{1}{d}}$.

\noi \hfill $\square$
\end{itemize}

\bigskip\noi
{\it Definition 1.6}

\noi Let $g = \sum\limits_{n = 0}^{\infty} a_n z^{\frac{n}{d}}$ be a
Puiseux series on ${\mathbf L}$. We call $d \in {\mathbb N}$ a
denominator of $g$. Note that $d$ is not unique.

\bigskip\noi {\it Definition 1.7}

\noi
\begin{itemize}
\item[a)] Let $r > 0$. We denote the set of all holomorphic
functions $B_{{\mathbf L}}(r) \to {\mathbf L}$ by ${\mathcal O}
{\mathbf L} (B_{{\mathbf L}}(r))$. We define an equivalence relation
$\equiv$ on $\bigcup\limits_{r > 0} {\mathcal O} {\mathbf L}
(B_{{\mathbf L}}(r))$ as follows: $f_1 \equiv f_2$ iff there is some
$r > 0$ such that $f_1, f_2 \in {\mathcal O} {\mathbf L}
(B_{{\mathbf L}}(r))$ and $f_1 \vert_{B_{{\mathbf L}}(r)} = f_2
\vert_{B_{{\mathbf L}}(r)}$. We let ${\mathcal O} {\mathbf L}_0$ be
the set of all $\equiv$ - equivalence classes.
\item[b)] We define ${\mathcal O} {\mathbf L}'_0 \subset {\mathcal O}
{\mathbf L}_0$ to be the set of all $\varphi \in {\mathcal O}
{\mathbf L}_0$ such that there is some $r > 0$, some $h \in
{\mathcal O} (B(0, r))$ with $h(0) = 0$ and $\vert h(z) \vert \leq
\frac{1}{2}$ for $\vert z \vert < r$, some $a \in {\mathbf L}$ and
some $k \in {\mathbb N}_0$ such that $\varphi (z) = {\mathbf m} (a,
{\mathbf m} ({\mathbf p}^k (z), 1 + h_{{\mathbf L}} (z)))$ for $z
\in B_{{\mathbf L}} (r)$. Note that $1 + h_{{\mathbf L}} (z) \in
{\mathbb C} \setminus {\mathbb R}_{\leq 0} \subset {\mathbf L}$ for
$\vert z \vert < r$. We write $r(\varphi)$, $h (\varphi)$, $a
(\varphi)$ and $k(\varphi)$ for the data above. We set ${\mathcal O}
{\mathbf L}^{\wedge}_0 := \{ \varphi \in {\mathcal O} {\mathbf L}'_0
\; \vert \; k (\varphi) \geq 1 \}$ and ${\mathcal O} {\mathbf L}^*_0
:= \{ \varphi \in {\mathcal O} {\mathbf L}'_0 \; \vert \; k
(\varphi) = 1 \}$. We define $s \colon {\mathcal O} {\mathbf L}'_0
\to {\mathcal O}_0$, $s (\varphi) (z) := a (\varphi)^1
z^{k(\varphi)} ( 1 + h (\varphi) (z))$.
\end{itemize}

\newpage\noi PROPOSITION 1.8

\begin{itemize}
\item[a)] {\it Let} $\varphi \in {\mathcal O} {\mathbf L}'_0$ {\it and} $\psi
\in {\mathcal O} {\mathbf L}^{\wedge}_0$. {\it Then} $\varphi \circ
\psi \in {\mathcal O} {\mathbf L}'_0$.
\item[b)]${\mathcal O} {\mathbf L}_0^*$ {\it is a group under composition.}
\end{itemize}

\bigskip\noi
{\it Proof}

\noi
\begin{itemize}
\item[a)] We see immediately that $a (\varphi \circ \psi) =
{\mathbf m} (a(\varphi)$, ${\mathbf p}^k (a(\psi)))$, that $k
(\varphi \circ \psi) = k (\varphi) k (\psi)$ and that
$$ h
(\varphi \circ \psi) = (1 + h (\psi))^{k (\varphi)} (1 + h (\varphi)
\circ s(\psi)) - 1.$$
\item[b)] Let $\varphi \in {\mathcal O} {\mathbf L}_0^*$. We show
that there is $\psi \in {\mathcal O} {\mathbf L}_0^*$ with $\varphi
\circ \psi = \id_{{\mathbf L}}$ and are done by a). We consider $f
:= s(\varphi)$. Then $f(0) = 0$ and $f'(0) = a (\varphi)^1 \ne 0$.
Hence $f$ is invertible at 0 and $f^{-1}$ has the form $f^{-1} (z) =
f'(0)^{-1} z (1 + \tilde{h} (z))$ with $\tilde{h} (0) = 0$. Let $b
\in {\mathbf L}$ with ${\mathbf m} (a (\varphi), b)) = (1, 0)$. Then
$\psi (z) := {\mathbf m} (b, {\mathbf m} (z, 1 + \tilde{h}_{{\mathbf
L}} (z)))$ fulfills the requirements.
\end{itemize}

\bigskip\noi PROPOSITION 1.9

\noi {\it Let} $\varphi, \psi \in {\mathcal O} {\mathbf L}_0^*$.

\noi a) $\vert \varphi (z) \vert \leq 2 \vert a (\varphi) \vert
\vert z \vert$ {\it and} $\vert \arg \varphi (z) - \arg  a (\varphi)
\vert \leq \vert \arg z \vert + \frac{\pi}{2}$
{\it for} $\vert z \vert < r (\varphi)$.\\
b) {\it We can choose} $r(\varphi \circ \psi) = \frac{1}{10}
\frac{\min \{ r(\varphi), r (\psi) \} }{\max \{ 1, \vert a (\psi)
\vert \} }$.

\bigskip\noi
{\it Proof}

\begin{itemize}
\item[a)] We have $\vert h (\varphi) (z) \vert < 1$ for $z \in B (0,
r(\varphi))$. Hence $\vert \varphi (z) \vert = \vert s (\varphi) (z)
\vert = \vert a (\varphi) \vert \; \vert z \vert \; \vert 1 + h
(\varphi) (z) \vert \leq 2 \vert a (\varphi) \vert \; \vert z
\vert$.

\noi We have $\vert \arg (1+ h (\varphi)_{{\mathbf L}} (z)) \vert
\leq \frac{\pi}{2}$ for $\vert z \vert < r (\varphi)$. This gives
the second part of a).
\item[b)] Let $r := r (\varphi)$, $s := r (\psi)$ and $t := \frac{1}{10} \frac{\min \{
r, s \} }{\max \{ 1, \vert a (\psi) \vert \} }$. By a) we have that
$s (\psi) (z) \in B(0, r)$ for $z \in B(0, t)$. Applying the maximum
principle we obtain that $\vert h(\varphi) (z) \vert \leq
\frac{\vert z \vert}{r}$ for $\vert z \vert < r$ and $\vert h (\psi)
(z) \vert \leq \frac{\vert z \vert}{s}$ for $\vert z \vert < s$.
With Proposition 1.8 we see that $h(\varphi \circ \psi) = h (\psi) +
h (\varphi) \circ s (\psi) + h (\psi) \cdot (h(\varphi) \circ s
(\psi))$ and obtain the claim by the above estimates and a).
\end{itemize}

\bigskip\noi PROPOSITION 1.10

\noi {\it Let $g$ be a Puiseux series convergent on $B_{{\mathbf
L}}(r) \subset {\mathbf L}$ with denominator $d$. Let $\varphi \in
{\mathcal O} {\mathbf L}^{\wedge}_0$. Then $g \circ \varphi$ is a
Puiseux series with denominator $d$ convergent on $B_{{\mathbf
L}}(s) \subset {\mathbf L}$ with} $s := \min \left\{ r (\varphi),
\left( \frac{r}{2 \vert a (\varphi) \vert}
\right)^{\frac{1}{k(\varphi)}} \right\}$.

\bigskip\noi
{\it Proof}

\noi As in Proposition 1.9 a) we see that $\vert \varphi (z) \vert
\leq 2 \vert a (\varphi) \vert \; \vert z \vert^{k(\varphi)}$ for
$\vert z \vert < r (\varphi)$. Hence $\varphi (z) \in B_{{\mathbf
L}}(r)$ for $\vert z \vert < s$ and therefore $g \circ \varphi \in
{\mathcal O} {\mathbf L} (B_{{\mathbf L}}(s))$. By binomial
expansion we obtain that $g \circ \varphi$ is a Puiseux series
convergent on $B_{{\mathbf L}}(t)$ for some $t \leq s$. We get the
claim by Proposition 1.5 a).

\bigskip\noi
{\it Definition 1.11} (compare with [27, p.6])

\noi Let $z$ be an indeterminate. A generalized log-power series in
$z$ is a formal expression $g(z) = \sum\limits_{\alpha \in {\mathbb
R}_{\geq 0}} a_{\alpha} P_{\alpha} (\log z) z^{\alpha}$ with
$a_{\alpha} \in {\mathbb C}$ and $P_{\alpha} \in {\mathbb C} [z]
\setminus \{ 0 \}$ monic such that $P_0 = 1$ and such that the
support of $g$ defined as $\supp (g) := \{ \alpha \in {\mathbb
R}_{\geq 0} \; \vert \; a_{\alpha} \ne 0 \}$ fulfils the following
condition: for all $R > 0$ the set $\supp (g) \cap [0, R]$ is
finite. We write ${\mathbb C} [[z^*]]^{\omega}_{\log}$ for the set
of generalized log-power series. These series are added and
multiplied by considering them as generalized power series with
logarithmic polynomials as coefficients. For $g \in {\mathbb C}
[[z^*]]^{\omega}_{\log}$ we set $\nu(g) := \min \supp (g)$. By
${\mathbb C} [[z^*]]^{\omega}$ we denote the subset of ${\mathbb C}
[[z^*]]^{\omega}_{\log}$ consisting of all $g \in {\mathbb C}
[[z^*]]^{\omega}_{\log}$ with $P_{\alpha} = 1$ for all $\alpha \in
{\mathbb R}_{\geq 0}$. By ${\mathbb C} [[z^*]]^{\omega,
\fin}_{\log}$ (resp. ${\mathbb C} [[z^*]]^{\omega, \fin}$) we denote
the set of all $g \in {\mathbb C} [[z^*]]^{\omega}_{\log}$ (resp.
${\mathbb C} [[z^*]]^{\omega}$) with finite support.

\bigskip\noi
{\it Convention}. From now on we omit the superscript $\omega$.

\bigskip\noi {\it Remark 1.12}

\begin{itemize}
\item[a)] The set ${\mathbb C} [[z^*]]_{\log}$ is a
${\mathbb C}$-algebra with ${\mathbb C} [[z^*]]$ as subalgebra.
\item[b)] Interpreting $\log z$ and $z^{\alpha}$ with Definition 1.1
we obtain that $g \in {\mathcal O} ({\mathbf L})$ for $g \in
{\mathbb C} [[z^*]]^{\fin}_{\log}$.
\end{itemize}

\newpage\noi {\it Definition 1.13} (compare with [27, Definition 2.3])

\noi A domain $W \subset {\mathbf L}$ of the Riemann surface of the
logarithm is a standard quadratic domain if there are constants $c,
C > 0$ such that
$$W = \left\{ (r, \varphi) \in {\mathbf L} \; \big\vert \; 0 < r < c \; \exp (-
C \sqrt{\vert \varphi \vert}) \right\}.$$

\noi A domain is called a quadratic domain if it contains a standard
quadratic domain.

\bigskip\noi
{\it Definition 1.14}

\noi Let $U \subset {\mathbf L}$ be a quadratic domain, let $f \in
{\mathcal O} (U)$ and let $g = \sum\limits_{\alpha \geq 0}
a_{\alpha} P_{\alpha} (\log z) z^{\alpha} \in {\mathbb C}
[[z^*]]_{\log}$. We say that $f$ has asymptotic expansion $g$ on $U$
and write $f \sim_U g$, if for each $R > 0$ there is a quadratic
domain $U_R \subset U$ such that
$$f(z) - \sum\limits_{\alpha \leq R} a_{\alpha} P_{\alpha} (\log z)
z^{\alpha} = o (\vert z \vert^R) \quad \mbox{as} \quad \vert z \vert
\longrightarrow 0 \quad \mbox{on} \quad U_R.$$

\noi We write $T f := g$. By ${\mathcal Q}^{\log} (U)$ we denote the
set of all $f \in {\mathcal O} (U)$ with an asymptotic expansion in
${\mathbb C} [[z^*]]_{\log}$. By ${\mathcal Q}(U)$ we denote the
subset of all $f \in {\mathcal Q}^{\log} (U)$ with $Tf \in {\mathbb
C} [[z^*]]$.

\bigskip\noi {\it Remark 1.15}

\noi If $f \in {\mathcal Q}^{\log}(U)$ for some quadratic domain $U$
then there is exactly one $g \in {\mathbb C} [[z^*]]_{\log}$ with $f
\sim_U g$, i.e. $Tf$ is well defined.

\bigskip\noi {\it Definition 1.16}

\noi We define an equivalence relation $\equiv$ on
$\bigcup\limits_{U \subset {\mathbf L} \; \mbox{\tiny quadr.}}
{\mathcal Q}^{\log} (U)$ as follows: $f_1 \equiv f_2$ if and only if
there is a quadratic domain $V \subset {\mathbf L}$ such that $f_1
\vert_V = f_2 \vert_V$. We let ${\mathcal Q}^{\log}$ be the set of
all $\equiv$ - equivalence classes. In the same way we obtain the
class ${\mathcal Q}$. Note that ${\mathcal Q} = {\mathcal Q}_1^1$ in
the notation of [27, Definition~5.1, Remarks~5.2 \& Definition 5.4].

\bigskip\noi {\it Remark 1.17}
\begin{itemize}
\item[a)] We will not distinguish between $f \in \bigcup\limits_{U \subset
{\mathbf L} \; \mbox{\tiny quadr.}} {\mathcal Q}^{\log} (U)$ and its
equivalence class in ${\mathcal Q}^{\log}$, which we also denote by
$f$. Thus ${\mathcal Q}^{\log} (U) \subset {\mathcal Q}^{\log}$
given a quadratic domain $U \subset {\mathbf L}$.
\item[b)] In the same way we define ${\mathcal Q} \subset {\mathcal Q}^{\log}$.
We have ${\mathcal Q}(U) \subset {\mathcal Q}$ for $U \subset
{\mathbf L}$ a quadratic domain.
\item[c)] Given a quadratic domain $U \subset {\mathbf L}$ the set ${\mathcal Q}^{\log}
(U)$ is a ${\mathbb C}$-algebra with ${\mathcal Q}(U)$ as a
subalgebra. Also, ${\mathcal Q}^{\log}$ is an algebra with
${\mathcal Q}$ as a subalgebra.
\item[d)] Given a quadratic domain $U \subset {\mathbf L}$ the well defined
maps $T \colon {\mathcal Q}^{\log} (U) \to {\mathbb C}
[[z^*]]_{\log}$, $f \mapsto Tf$, and $T \colon {\mathcal Q}(U) \to
{\mathbb C} [[z^*]]$, $f \mapsto Tf$, are homomorphisms of ${\mathbb
C}$-algebras. Also the induced maps $T \colon {\mathcal Q}^{\log}
\to {\mathbb C} [[z^*]]_{\log}$, $f \mapsto Tf$, and $T \colon
{\mathcal Q} \to {\mathbb C} [[z^*]]$, $f \mapsto Tf$, are
homomorphisms of ${\mathbb C}$-algebras.
\end{itemize}

\bigskip\noi
PROPOSITION 1.18

\noi {\it Let} $U \subset {\mathbf L}$ {\it be a quadratic domain.
The homomorphism} $T \colon {\mathcal Q}^{\log} (U) \to {\mathbb C}
[[z^*]]_{\log}$ {\it is injective. Therefore the homomorphism} $T
\colon {\mathcal Q}^{\log} \to {\mathbb C} [[z^*]]_{\log}$ {\it is
injective.}

\bigskip\noi
{\it Proof}

\noi See Ilyashenko [22, Theorem 2 p.23] and [27, Proposition 2.8].
\hfill $\square$

\bigskip\noi
PROPOSITION 1.19

\begin{itemize}
\item[a)] {\it Let $g(z) = \sum\limits_{n = 0}^{\infty} a_n z^{\frac{n}{d}}$
be a Puiseux series convergent on $B_{{\mathbf L}}(r) \subset
{\mathbf L}$. Then $g \in {\mathcal Q} (B_{{\mathbf L}}(r))$ with}
$Tg = \sum\limits_{n = 0}^{\infty} a_n z^{\frac{n}{d}}$.
\item[b)] {\it Let} $f \in {\mathcal Q}^{\log}$ {\it and let} $\rho > 0$. {\it Then} $f
\circ {\mathbf p}^{\rho} \in {\mathcal Q}^{\log}$. {\it If} $f \in
{\mathcal Q}$ {\it then} $f \circ {\mathbf p}^{\rho} \in {\mathcal
Q}$.
\item[c)] {\it Let} $f \in {\mathcal Q}^{\log}$ {\it and let} $\psi \in
{\mathcal O} {\mathbf L}^{\wedge}_0$. {\it Then} $f \circ \psi \in
{\mathcal Q}^{\log}$. {\it If} $f \in {\mathcal Q}$ {\it then} $f
\circ \psi \in {\mathcal Q}$.
\end{itemize}

\newpage\noi {\it Proof}

\noi
\begin{itemize}
\item[a)] is a consequence of Proposition 1.5 b).
\item[b)] Let $c, C > 0$ such that $f \in {\mathcal Q}^{\log} (U)$
where $U := \{ (r, \varphi) \in {\mathbf L} \; \vert \; 0 < r < c\,
\exp (-C \sqrt{\vert \varphi \vert}) \}$. Let $d :=
c^{\frac{1}{\rho}}$, $D := \frac{C}{\sqrt{\rho}}$ and
$$W := \left\{
(r, \varphi) \in {\mathbf L} \; \big\vert \; 0 < r < d \; \exp (- D
\sqrt{\vert \varphi \vert}) \right\}.$$

\noi Then we see that ${\mathbf p}^{\rho} (W) \subset U$ and
therefore $f \circ \varphi \in {\mathcal O} (W)$. The claim follows
from the following observation: let $\alpha \geq 0$ and $m \in
{\mathbb N}_0$. Then $(z^{\alpha} (\log z)^m) \circ {\mathbf
p}^{\rho} = \rho^m z^{\alpha \rho} (\log z)^m$.
\item[c)] Let $k := k (\psi)$, $b := {\mathbf
p}^{\frac{1}{k}} (a(\psi))$ and $\tilde{h} (z) := \root k \of {1 + h
(\psi) (z)} -1$. Then $\chi := {\mathbf m} (b, {\mathbf m} (z, 1 +
\tilde{h}_{{\mathbf L}} (z)) \in {\mathcal O} {\mathbf L}_0^*$ and
$\psi = {\mathbf p}^k \circ \chi$. Hence we can assume by b) that
$\psi \in {\mathcal O} {\mathbf L}_0^*$.

\noi Let $c, C > 0$ such that $f \in {\mathcal Q}^{\log} (U)$ where
$$U :=  \left\{ (r, \varphi) \in {\mathbf L} \; \big\vert \; 0 < r <
c \; \exp (- C \sqrt{\vert \varphi \vert}) \right\}.$$

\noi Let $r_0 := r (\psi)$ and $z = (r, \varphi) \in B_{{\mathbf L}}
(r_0)$. Then $\vert \psi (z) \vert \leq 2 \vert a (\psi) \vert r$
and $\vert \arg \psi (z) \vert \leq \vert \arg a (\psi) \vert +
\vert \varphi \vert + \frac{\pi}{2}$ by Proposition 1.9 a). Let
$$d := \min \left\{ \frac{c}{2 \vert a (\psi) \vert} \; \exp (- C
\sqrt{\vert \arg a (\psi) \vert + \frac{\pi}{2}}), r_0 \right\}$$

\noi and $W \colon = \left\{ (r, \varphi) \in {\mathbf L} \; \vert
\; 0 < r < d \; \exp (- C \sqrt{\vert \varphi \vert}) \right\}
\subset B_{{\mathbf L}}(r_0)$. Then $\psi (W) \subset U$ and
therefore $f \circ \psi \in {\mathcal O} (W)$. Let $Tf =:
\sum\limits_{\alpha \geq 0} a_{\alpha} P_{\alpha} (\log z)
z^{\alpha}$. Let $R > 0$. Then there is a quadratic domain $U_R$
such that
$$f(z) - \sum\limits_{\alpha \leq R} a_{\alpha} P_{\alpha} (\log z)
z^{\alpha} = o (\vert z \vert^R) \quad \mbox{as} \quad \vert z \vert
\longrightarrow 0 \quad \mbox{on} \quad U_R.$$

\noi As above we find some quadratic domain $W_R \subset B_{{\mathbf
L}}(r_0)$ such that $\psi (W_R) \subset U_R$. With Proposition 1.9
a) we obtain
$$f \circ \psi (z) - \left( \sum\limits_{\alpha \leq R} a_{\alpha}
P_{\alpha} (\log z) z^{\alpha} \right) \circ \psi = o (\vert z
\vert^R) \quad \mbox{as} \quad \vert z \vert \longrightarrow 0 \quad
\mbox{on} \quad W_R.$$

\noi Hence the claim follows from the next lemma. The second part of
it will be used in section 3 below.
\end{itemize}

\bigskip\noindent LEMMA 1.20

\noindent Let $\alpha > 0$ and let $m \in {\mathbb N}_0$. Let
$\varphi \in {\mathcal O} L_0^{\wedge}$. Then there are power series
$g_1, \dots, g_m$ convergent on $B_{{\mathbf L}} (r(\varphi))$ such
that $(z^{\alpha} (\log z)^m) \circ \varphi = z^{k (\varphi) \alpha}
\sum\limits_{\ell = 0}^m g_{\ell} (z) (\log z)^{\ell}$. If $\varphi
\in {\mathcal O} L_0^*$ with $\vert a (\varphi) \vert = 1$ then
$\vert g_{\ell} (z) \vert \leq 2^{m + \alpha} (\vert \arg a
(\varphi) \vert + 3)^m$ for every $\ell = 1, \dots, m$.

\bigskip\noi
{\it Proof}

\noindent Let $a := a (\varphi)$, $k := k (\varphi)$ and $h := h
(\varphi)$. By the binomial formula we obtain
$$(z^{\alpha} (\log z)^m) \circ \varphi = a^{\alpha} z^{k \alpha}
(1 + h(z))^{\alpha} \sum\limits_{\ell = 0}^m {m \choose \ell} (\log
{\mathbf m} (a, 1 + h_{{\mathbf L}} (z)))^{m - \ell} (\log {\mathbf
p}^k (z))^{\ell}.$$

\noindent Since $\log {\mathbf p}^k (z) = k \log z$ we see that
$$g_{\ell} := k^{\ell} {m \choose \ell} a^{\alpha} (1 + h
(z))^{\alpha} (\log {\mathbf m} (a, 1 + h_{{\mathbf L}}
(z)))^{m-\ell} \leqno(\ast)$$

\noindent fulfils the claim.

\medskip\noi
Assume that $\varphi \in {\mathcal O} L_0^*$ with $\vert a (\varphi)
\vert = \vert a \vert = 1$. We have $k = 1$ and $\vert a^{\alpha}
\vert = 1$. Since $\vert h_{{\mathbf L}} (z) \vert \leq \frac{1}{2}$
for $z \in B_{{\mathbf L}} (r(\varphi))$ we obtain $\vert (1 +
h_{{\mathbf L}} (z))^{\alpha} \vert \leq \left( \frac{3}{2}
\right)^{\alpha} \leq 2^{\alpha}$. Since ${m \choose \ell} \leq 2^m$
for all $\ell$ the non-logarithmic terms in ($\ast$) are bounded
from above by $2^{m+\alpha}$. We have
$$\begin{array}{l}
\log {\mathbf m} (a, 1 + h_{{\mathbf L}} (z)) = \\
\qquad = \log a + \log (1 + h (z)) \\
\qquad = i \arg a + \log \vert 1 + h (z) \vert + i \arg (1 + h
(z)). \\
\end{array}$$

\noindent Since $\frac{1}{2} \leq \vert 1 + h (z) \vert \leq
\frac{3}{2}$ for $z \in B_{{\mathbf L}} (r(\varphi))$ we see that
$\vert \log \vert 1 + h (z)) \vert \leq \log 2 \leq 1$ for $z \in
B_{{\mathbf L}} (r(\varphi))$. Since $\vert \arg (1 + h (z)) \vert
\leq \frac{\pi}{2}$ we conclude finally that $\vert \log {\mathbf m}
(a, 1 + k (z)) \vert \leq \arg a + 1 + \frac{\pi}{2} \leq \arg a +
3$. This gives the claim. \hfill $\square$

\newpage\noindent {\bf 2. Angles and domains with an analytic corner}

\medskip\noindent
We introduce the notion of an angle. This allows us to formulate
Theorem A and B in a precise way.

\noi The proof of Theorem A will be reduced in Section 3 to domains
with an analytic corner as considered by Wasow (see [39]). We give
the definition.

\noi Sets definable in o-minimal structures are subsets of cartesian
products of the reals . So for Theorem~B we have to stay on the
reals. For Theorem~A we have to go out on the Riemann surface of the
logarithm. We have to distinguish the ambient spaces carefully and
we describe the lifting process. Finally we extend the result of
Wasow on the existence of asymptotic expansion (see [39]) to the
case where the boundary is given by Puiseux series.

\bigskip\noi {\it Remark 2.1}

\noi Let $\Omega$ be a bounded and subanalytic domain in ${\mathbb
R}^n$. Let $x \in \partial \Omega := \overline{\Omega} \setminus
\Omega$. Then the germ of $\Omega$ at $x$ has finitely many
connected components. More precisely we have the following: there is
$k \in {\mathbb N}_0$ such that for all neighbourhoods $V$ of $x$
exactly $k$ of the components of the set $\Omega \cap V$ have $x$ as
a boundary point.

\bigskip\noi
{\it Remark 2.2}

\noi Let $\Omega \subset {\mathbb R}^2$ be a bounded and
semianalytic domain without isolated boundary points. Let $x \in
\partial \Omega$ and let $C$ be a connected component of the germ of
$\Omega$ at $x$. Then the germ of the boundary of $\partial C$ at
$x$ is given by (the germs of) one or two semianalytic curves (see
[4, Theorem 6.13]). So the interior angle of $C$ at $x$, denoted by
$\sphericalangle_{x} C$, is well defined; it takes value in $[0, 2
\pi]$. If the germ of $\Omega$ at $x$ is connected we write
$\sphericalangle_{x} \Omega$.

\bigskip\noi
{\it Definition 2.3}

\noi Let $\Omega \subset {\mathbb R}^2$ be a bounded semianalytic
domain without isolated boundary points.
\begin{itemize}
\item[a)] A point $x \in \partial \Omega$ is a {\it singular boundary point}
if $\partial \Omega$ is not a real analytic manifold at $x$.
\item[b)] We set $Sing (\partial \Omega) := \{ x \in \partial \Omega \;
\vert \; x$ is a singular boundary point of $\partial \Omega \}$.
\item[c)] Let $x \in \partial \Omega$. We set $\sphericalangle (\Omega, x) := \{
\sphericalangle_{x} C \; \vert \; C$ is a component of the germ of
$\Omega$ at $x$ and $x \in \Sing (\partial C) \}$.
\end{itemize}

\bigskip\noi {\it Remark 2.4}

\noi Let $\Omega \subset {\mathbb R}^2$ be a bounded semianalytic
domain without isolated boundary points.
\begin{itemize}
\item[a)] $\Sing (\partial \Omega)$ is finite by analytic cell
decomposition (see [13, pp.508-509]).
\item[b)] Let $x \in \partial \Omega$. Then $\sphericalangle (\Omega, x) = \emptyset$
if and only if $x \not\in \Sing (\partial C)$ for all components $C$
of the germ of $\Omega$ at $x$. This is especially the case if $x
\not\in \Sing (\partial \Omega)$.
\end{itemize}

\bigskip\noi {\it Definition 2.5} (compare with the introduction of [39])

\noi \begin{itemize}
\item[a)] We say that a domain $D \subset
{\mathbb C}$ with $0 \in
\partial D$ has an analytic corner (at 0) if the boundary of $D$ at
0 is given by two analytic curves which are regular at 0 and if $D$
has an interior angle greater than 0 at $0$.
More precisely, the following holds.\\
There are holomorphic functions $\psi, \chi \in {\mathcal O}
(B(0,1))$ with $\psi (0) = \chi (0) = 0$ and $\psi' (0) \chi' (0)
\ne 0$ such that for $\Gamma := \psi ([0, 1[)$ and $\Gamma := \chi
([0, 1[)$ the following holds:

\medskip\noi
i) $\partial D \cap B (0, r) = (\Gamma \cup \Gamma') \cap B (0, r)$
for
some $r > 0$.\\
ii) The interior angle $\sphericalangle D \in [0, 2 \pi ]$ of
$\partial D$ at 0 is greater than 0.

\medskip\noi
Note that possibly $\Gamma_1 = \Gamma_2$ if $\sphericalangle D = 2
\pi$. Otherwise we may assume that $\Gamma_1 \cap \Gamma_2 = \{ 0
\}$. We call $\psi$ and $\chi$ $D$-describing functions. Note that
$D \cap B (0, r)$ is semianalytic for all sufficiently small $r
> 0$.
\item[b)] Let $h \colon \partial D \to {\mathbb R}$ be a continuous
boundary function. We call $h$ a corner function (at 0) if the
following holds: There is some $\varepsilon$ with $0 < \varepsilon
\leq 1$ and there are real Puiseux series $g_{0}$, $g_{1} \colon ]
0, \varepsilon [ \longrightarrow {\mathbb R}$ such that $h \circ
\psi (t) = g_{0} (t)$, $h \circ \chi (t) = g_{1} (t)$ for $0 \leq t
< \varepsilon$ (note that $h_0 := g_0 \circ \psi^{-1}$ and $h_1 :=
g_1 \circ \chi^{-1}$ are Puiseux series convergent on $B_{{\mathbf
L}} (s)$ (compare with Proposition 1.10) and that $h = h_0$ on
$\Gamma \cap B(s)$ and $h = h_1$ on $\Gamma' \cap B(s)$ for some $s
> 0$).
\end{itemize}

\bigskip\noi {\it Definition 2.6}

\noi
\begin{itemize}
\item[a)] We set $e:= \exp \circ \log \colon {\mathbf L} \to
{\mathbb C}^*$, $(r, \varphi) \mapsto r e^{i \varphi}$. Via the
identification of ${\mathbb C} \setminus {\mathbb R}_{\leq 0}$ with
${\mathbb R}_{> 0} \times ] - \pi, \pi[$ (see Definition~1.1) we see
that $e \vert_{{\mathbb C} \setminus {\mathbb R}_{\leq 0}} =
\id_{{\mathbb C} \setminus {\mathbb R}_{\leq 0}}$.
\item[b)] Let $A \subset {\mathbb C}^*$ be a set. We say that $A$
can be embedded in ${\mathbf L}$ if there is a set $B \subset
{\mathbf L}$ such that $e \vert_B$ is injective and $e(B) = A$. We
call $B$ an embedding of $A$ in ${\mathbf L}$.
\item[c)] Let $A \subset {\mathbb C}^*$ be embeddable in
${\mathbf L}$ and let $f \colon A \to {\mathbb C}$ be a function.
Let $B$ be an embedding of $A$. We set $f_B \colon = f \circ e
\colon B \to {\mathbb C}$.
\end{itemize}

\bigskip\noi
{\it Remark 2.7}

\noi
\begin{itemize}
\item[a)] Let $U \subset {\mathbb C}^*$ be an embeddable domain.
Then an embedding of $U$ is a domain in ${\mathbf L}$.
\item[b)] Let $\Omega$ be a semianalytic domain such that $0$ is a
non-isolated boundary point of $\Omega$. Then every component of the
germ of $\Omega$ at 0 has a semianalytic representative which is
embeddable in ${\mathbf L}$.
\item[c)] Let $A \subset {\mathbb C}^*$ be embeddable in
${\mathbf L}$ and let $f \colon A \to {\mathbb C}$ be a function.
Let $B_1, B_2$ be embeddings of $A$. Then there is some $k \in
{\mathbb Z}$ such that $B_2 = \left\{ {\mathbf m} (a, z) \; \vert \;
z \in B_1 \right\}$ and $f_{B_2} (z) = f_{B_1} ({\mathbf m} (a^{-1},
z))$ where $a := (1, 2 k \pi)$ and $a^{-1} := (1, - 2 k\pi)$.
\end{itemize}

\bigskip\noi
{\it Definition 2.8}

\noi Let $U \subset {\mathbb C}^*$ be an embeddable domain and let
$f \colon U \to {\mathbb C}$ be a function. Let $W \subset {\mathbf
L}$ be a domain and $F \colon W \to {\mathbb C}$ be a function. We
say that $F$ extends $f$ if there is an embedding $V$ of $U$ in
${\mathbf L}$ with $V \subset W$ such that $F$ extends $f_W$.

\bigskip\noi THEOREM 2.9

\noi {\it Let $D$ be a bounded and simply connected domain with an
analytic corner at $0$. Let $h \colon \partial D \to {\mathbb R}$ be
a continuous function which is a corner function at $0$. Let $u$ be
the Dirichlet solution for $h$. Let $r > 0$ such that $D \cap B (0,
r)$ is embeddable and let $D'$ be an embedding of $D \cap B (0, r)$.
Then there is some $f = f (D') \in {\mathcal O} (D') \cap C
(\overline{D'})$ with} $\Re f = u_{D'}$ {\it and some} $g = g (D') =
\sum\limits_{\alpha \geq 0} a_{\alpha} P_{\alpha} (\log z)
z^{\alpha} \in {\mathbb C} [[z^*]]_{\log}$ {\it such that for} $R
> 0$
$$f(z) - \sum\limits_{\alpha \leq R} a_{\alpha} P_{\alpha} (\log z)
z^{\alpha} = o (\vert z \vert^R) \quad as \quad \vert z \vert \to 0
\quad on \quad \overline{D'}.$$

\noi {\it If} $\sphericalangle D / \pi \in {\mathbb R} \setminus
{\mathbb Q}$ {\it then} $g \in {\mathbb C} [[z^*]]$.

\bigskip\noi
{\it Proof}

\noi If the boundary functions $g_0$ and $g_1$ are power series we
work exactly in the setting of [39] and we obtain the statement by
[39, Theorem~3 \& Theorem~4]. Looking carefully at the proofs of
these theorems we can generalize the results to Puiseux series: as
in [32] we restrict ourselves to the case where $\alpha :=
\sphericalangle D/\pi$ is irrational (compare with [39, p.55]) and
use the notation introduced there. Up to the end of [39, p.53] the
arguments are literally the same. But in our case the function
$\phi_2 (s)$ is in general a convergent real Puiseux series $\phi_2
(s) = \phi (s) = \sum\limits_{n = n_0}^{\infty} a_n s^{\frac{n}{d}}$
where $n_0 > 0$ and $a_{n_0} \ne 0$. By (4.7) of [39] we have to
consider $\frac{\partial \mu (\xi, 0)}{\partial \xi}$ which is given
by (4.13) of [39] as
$$\frac{\partial \mu
(\xi, 0)}{\partial \xi} = \phi' (s) / \frac{d \xi}{ds}$$

\noi By (4.12) of [39] we have $s = \psi (\xi) = \xi^{\alpha}
K^{\ast\ast\ast} (\xi)$ where $K^{\ast\ast\ast} (\xi) \sim
\sum\limits_{k, \ell \geq 0} a_{k \ell} \xi^{k + \ell \alpha}$ on
the real line and where $a_{00} \ne 0$. Hence by binomial expansion
$$\phi' (\psi(\xi)) \sim \xi^{(\frac{n_0}{d}-1)\alpha}
\sum\limits_{k, \ell \geq 0} b_{k \ell} \xi^{k + \frac{\ell}{d}
\alpha}$$

\noi where $b_{00} \ne 0$. Since $\frac{ds}{d \xi} \sim
\xi^{\alpha-1} \sum\limits_{k, \ell \geq 0} c_{k \ell} \xi^{k + \ell
\alpha}$ where $c_{00} \ne 0$ (compare with (4.15) of [39]) we see
that
$$\frac{\partial \mu (\xi, 0)}{\partial \xi} \sim
\xi^{\frac{n_0}{d} \alpha-1} \sum\limits_{k, \ell \geq 0} d_{k \ell}
\xi^{k + \frac{\ell}{d} \alpha}.$$

\noi Note that $\frac{n_0}{d} \alpha - 1 > - 1$ and that all
exponents are irrational. Therefore Lemma 1, 2 \& 3 of [39] can be
applied to formula (4.7) of [39] and we can finish the proof by the
arguments of [39, pp.54-55]. \hfill $\square$

\newpage\noi {\bf 3. Proofs of Theorem A and B}

\medskip\noindent
We prove Theorem~A by reducing the problems to the case of a domain
with an analytic corner at 0. Theorem~B is then deduced from
Theorem~A by applying the results and methods of [27].

\bigskip\noi
{\it Definition 3.1}

\noi Let $\Omega \subset {\mathbb R}^2$ be a bounded and
semianalytic domain without isolated boundary points such that $0
\in \partial  \Omega$. We call $C \subset \Omega$ a corner component
of $\Omega$ at $0$ if $C$ is a semianalytic, simply connected and
embeddable representative of a connected component of the germ of
$\Omega$ at 0 such that the germ of $\partial C$ at 0 consists of
one or two semianalytic curves.

\medskip\noi
Theorem A gets now the following precise form:

\bigskip\noi THEOREM 3.2

\noi {\it Let} $\Omega \subset {\mathbb R}^2$ {\it be a bounded and
semianalytic domain without isolated boundary points with} $0 \in
\partial \Omega$. {\it Let} $h \colon \partial \Omega \to {\mathbb R}$ {\it be
a continuous and semianalytic boundary function and let $u$ be the
Dirichlet solution for $h$. Let $C$ be a corner component of
$\Omega$ at $0$. If $\sphericalangle_0 C > 0$ then there is a
quadratic domain $U \subset {\mathbf L}$ and an} $f \in {\mathcal
Q}^{\log} (U)$ {\it such that} $\Re f$ {\it extends $u \vert_C$. If
$\sphericalangle_0 C /\pi \in {\mathbb R} \setminus {\mathbb Q}$
then} $f \in {\mathcal Q} (U)$.

\medskip\noi
Note that Theorem 3.2 does not depend on the chosen embedding by
Definition 1.13, Definition 1.14 and Remark 2.7 c). Theorem 3.2 can
be deduced from the following:

\bigskip\noi
THEOREM 3.3

\noi {\it Let $D$ be a bounded, simply connected and embeddable
domain with an analytic corner at $0$. Let $h \colon \partial D \to
{\mathbb R}$ be a continuous function which is a corner function at
$0$. Let $u$ be the Dirichlet solution for $h$. Then there is a
quadratic domain $U \subset {\mathbf L}$ and an} $f \in {\mathcal
Q}^{\log} (U)$ {\it such that} $\Re f$ {\it extends $u$. If
$\sphericalangle D / \pi \in {\mathbb R} \setminus {\mathbb Q}$
then} $f \in {\mathcal Q} (U)$.

\bigskip\noi
{\it Proof of Theorem 3.2 from Theorem 3.3}

\noi We can assume that $\partial C \cap \partial \Omega$ consists
of two semianalytic branches $\Gamma$ and $\Gamma'$. Moreover, we
may assume that after some rotation there is a convergent Puiseux
series $g \colon [0, \delta [ \to {\mathbb R}$ such that $\Gamma =
\left\{ (t, g (t)) \; \vert \; 0 \leq t \leq \gamma \right\}$ for
some $0 < \gamma < \delta$. This can be achieved by analytic cell
decomposition (see [13, pp.508-509]) and the fact that bounded
semianalytic functions in one variable are given by Puiseux series
(see [9, p.192]). There is some $d \in {\mathbb N}$ and some
convergent real power series $f \colon ] - \delta^d, \delta^d [ \to
{\mathbb R}$ such that $g(t) = f (t^{\frac{1}{d}})$. Hence $\Gamma =
\left\{ (t^d, f (t)) \; \vert \; 0 \leq t \leq \gamma^d \right\}$.
We consider $\psi \colon B (0, \delta^d) \to {\mathbb C}$, $z
\mapsto z^d + i f (z)$. Then $\psi \in {\mathcal O} (B (0,
\delta^d))$ and $\Gamma = \psi ([0, \gamma^d])$. Arguing similarly
for $\Gamma'$ we find (after back-rotation and some dilation)
holomorphic functions $\psi, \chi \in {\mathcal O} (B(0, 1))$ with
$\psi (0) = \chi (0)$ such that $\Gamma = \psi ([ 0, 1[)$ and
$\Gamma' = \chi ([ 0, 1'[)$. Let $m, n \in {\mathbb N}$ such that
$\psi (z) = z^{m} \widehat{\psi} (z)$, $\chi (z) = z^{n}
\widehat{\chi} (z)$ with $\widehat{\psi}$, $\widehat{\chi} \in
{\mathcal O} (B(0, 1))$ and $\widehat{\psi} (0) \widehat{\chi} (0)
\ne 0$. We can replace $\chi (z)$ by $\chi (z^{m})$ and can
therefore assume that $m$ divides $n$.\\
We choose an embedding of $C$ in ${\mathbf L}$ which we denote again
by $C$. By the above we have functions from ${\mathcal O} {\mathbf
L}_0^{\wedge}$, denoted again by $\psi$ and $\chi$, which are
defined on $B_{{\mathbf L}}(1)$ such that $k (\psi) = m, k (\chi) =
n$ and $\partial C \cap B_{{\mathbf L}}(r) = \Gamma \cup \Gamma'$
for some $r > 0$ where $\Gamma := \psi ([0, 1[)$ and $\Gamma' :=
\chi ([ 0, 1'[)$. We apply finitely many elementary manipulations
from Section 1 to $u$ and $C$. We obtain functions $u_i$ and domains
$C_i$, $i = 1, 2, 3$, such that $C_3$ is an embedding of a domain
with an analytic corner. This enables us to apply Theorem 3.3. The
domains $C_i$ allow a similar description as $C$. We denote the data
describing $C_i$ by $\Gamma_i, \Gamma'_i, \psi_i, \chi_i$.

\begin{itemize}
\item[1)] We consider $u_1 \colon C_1 \to {\mathbb
R}$, $z \mapsto u ({\mathbf p}^{k (\psi)} (z))$, where $C_1 :=
{\mathbf p}^{\frac{1}{k (\psi)}} (C)$. Since $k (\psi)$ divides $k
(\chi)$ we see that $C_1$ has similar properties as $C$ but
additionally $k (\psi_1) = 1$.
\item[2)] We may choose a priori $C$ and $r(\psi_1^{-1})$
such that $\psi_1^{-1}$ is injective on $B_{{\mathbf
L}}(r(\psi_1^{-1}))$ and $C_1 \subset B_{{\mathbf L}} (r
(\psi_1^{-1}))$ (compare with Definition 1.7 b) and Proposition 1.8
b)).\\
We consider $u_2 \colon C_2 \to {\mathbb R}$, $z \mapsto u_1 (\psi_1
(z))$, where $C_2 := \psi_1^{-1} (C_1 )$. Then $C_2$ has similar
properties as $C_1$ but additionally $\Gamma_2 \subset {\mathbb
R}_{> 0} \times \{ \varphi \}$ for some $\varphi \in {\mathbb R}$.
\item[3)] We consider $u_3 \colon C_3 \to
{\mathbb R}$, $z \mapsto u_2 ({\mathbf p}^{k (\chi_2)})$, where $C_3
:= {\mathbf p}^{\frac{1}{k (\chi_2)}} (C_2)$. Then additionally $k
(\chi_3) = 1$.
\end{itemize}

\medskip\noi
By construction $e \vert_{C_3}$ is injective and $e(C_3)$ is a
bounded, simply connected and embeddable domain with an analytic
corner at $0$. Moreover $u_e := u_3 \circ (e \vert_{C_3})^{-1}$ is a
harmonic function on $e (C_3)$ which has a continuous extension to
$\overline{e (C_3)}$. With $h_e$ we denote the extension to the
boundary of $e (C_3)$. With Proposition 1.8 and
Proposition 1.10 we see that $h_e$ is a corner function at $0$.\\
By Theorem 3.3 there is a quadratic domain $U$ and a $g \in
{\mathcal Q}^{\log} (U)$ such that $\Re\, g$ extends $u_3$. If
$\sphericalangle e (C_3)/ \pi \in {\mathbb R} \setminus {\mathbb Q}$
then $g \in {\mathcal Q} (U)$. Since $\sphericalangle e (C_3) = (k
(\psi) k (\chi_2))^{-1} \sphericalangle_0 C$ we get that $g \in
{\mathcal Q} (U)$ if $\sphericalangle_0 C / \pi \in {\mathbb R}
\setminus {\mathbb Q}$. By construction $u = u_3 \circ {\mathbf
p}^{\frac{1}{k (\chi_2)}} \circ \psi_1^{-1} \circ {\mathbf
p}^{\frac{1}{k (\psi)}}$.\\
By Proposition 1.19 we obtain that $f := g \circ {\mathbf
p}^{\frac{1}{k (\chi_2)}} \circ \psi_1^{-1} \circ {\mathbf
p}^{\frac{1}{k (\psi)}} \in {\mathcal Q}^{\log}$ and $f \in
{\mathcal Q}$ if $\sphericalangle_0 C / \pi \in {\mathbb R}
\setminus {\mathbb Q}$. This $f$ fulfils the requirements. \hfill
$\square$

\bigskip\noi
It remains to prove Theorem 3.3. We prove it by doing reflections at
analytic curves infinitely often. We use the fact that the given
boundary function is defined at 0 by convergent Puiseux series which
extend to the Riemann surface of the logarithm (see Remark 1.4). To
motivate the technical statements of the upcoming proofs we give the
following example for the reflection principle involved.

\bigskip\noi {\it Example 3.4}

\noi Let $R > 0$ and let $\chi \colon B(0, R) \to {\mathbb C}$ be an
injective holomorphic function with $\chi (0) = 0$, $\chi' (0) > 0$
and $\Gamma := \chi (] 0, R[) \subset {\mathbb C}_+ := \left\{ z \in
{\mathbb C} \; \vert \; \Re z > 0 \right\}$ (note that $\Gamma$ is
tangent to ${\mathbb R}_{> 0}$ at 0). Let $R' > 0$ such that $\chi
(B(0, R)) \supset B(0, R')$. Then $B(0, R')\!\setminus\!\Gamma \cap
{\mathbb C}_+$ has two components, let $V$ be one of them. Let $f
\colon V \to {\mathbb C}$ be a holomorphic function which has a
continuous extension to $\Gamma \cap B(0, R')$ such that the
following holds: there is a convergent Puiseux series $h \colon B(0,
R') \setminus {\mathbb R}_{\leq 0} \to {\mathbb R}$ such that $\Re f
= h$ on $\Gamma \cap B(0, R')$. Then there is some $R''$ with $0 <
R'' < R'$ such that $f$ has a holomorphic extension to $B(0, R'')$
given by
$$z \longmapsto \left\{ \begin{array}{lll}
f(z) & & z \in V, \\
& \mbox{if} &     \\
\overline{(f-h) \circ \chi \circ \overline{\chi^{-1} (z)} + h (z)} &
& z \not\in V\\
\end{array} \right.$$

\noi (here $\overline{z}$ denotes the complex conjugate of a complex
number $z$). This can be seen in the following way: $f_1 := f-h \in
{\mathcal O} (V)$ has a continuous extension to $\Gamma \cap B(0,
R')$ with $\Re f_1 = 0$ on $\Gamma \cap B(0, R')$. Then $f_2 := f_1
\circ \chi \in {\mathcal O} (W)$ with $W := \chi^{-1} (V)$ has a
continuous extension to $I := \chi^{-1} (\Gamma \cap B(0, R'))
\subset {\mathbb R}$ with vanishing real part there. Therefore $f_3
\colon W \cup I \cup W^* \to {\mathbb C}$ defined by
$$z \longmapsto \left\{ \begin{array}{lll}
f_2 (z) & & z \in W \cup I, \\
& \mbox{if} &               \\
- \overline{f_2 (\overline{z})} & & z \in W^*, \\
\end{array} \right. $$

\noi where $W^* := \{ \overline{z} \; \vert \; z \in W \}$ is a
holomorphic extension of $f_2$. Then there is some $R''$ with $0 <
R'' \leq R'$ such that $f_4 := f_3 \circ \chi^{-1} \in {\mathcal O}
(B(0, R''))$ extends $f_2$. So $f_5 : = f_4 + h \in {\mathcal O}
(B(0, R''))$ extends $f$.

\bigskip\noi {\it Remark 3.5}

\noi We define the conjugate $\tau \colon {\mathbf L} \to {\mathbf
L}$ by $\tau (r, \varphi) := (r, - \varphi)$. We obtain immediately
the following: let $g$ be a Puiseux series convergent on
$B_{{\mathbf L}}(r)$ with denominator $d$. Then $\overline{g \circ
\tau}$ is a Puiseux series convergent on $B_{{\mathbf L}}(r)$ with
denominator $d$. Let $\varphi \in {\mathcal O} {\mathbf L}'_0$. Then
$\psi := \tau \circ \varphi \circ \tau \in {\mathcal O} {\mathbf
L}'_0$ with $r (\psi) = r (\varphi)$, $k (\psi) = k(\varphi)$ and
$\vert a (\psi) \vert = \vert a (\varphi) \vert$.

\medskip\noi
The proof of Theorem 3.3 relies on the iteration of the reflection
process from Example 3.4. Its description involves inversion,
conjugation and iteration of certain holomorphic functions and
summation of Puiseux series (see ($\dagger$) in the proof below). We
separate the proof into two steps. In the first step we show the
extension of the Dirichlet solution to a quadratic domain (compare
with [21, Proposition 3] and the subsequent remarks there). In the
second step we show the extension of the asymptotic development. We
use the classes of functions from Section~1 to establish the
dynamical system ($\dagger$) and we use the estimates for these
classes from Section~1 to control it to get the desired properties.

\bigskip\noi {\it Proof of Theorem 3.3}

\noi We show the claim in two steps:

\medskip\noi{\bf Step 1:} There is a quadratic domain $U \subset
{\mathbf L}$ and an $f \in {\mathcal O} (U)$ such that $\Re f$
extends $u$.

\medskip\noi
Proof of Step 1:\\
We choose an embedding of $D$ which we also denote by $D$. We also
write $u$ for $u_D$ and $h$ for $u \vert_{\partial D}$. Let $\psi,
\chi$ be $D$-describing (see Definition 2.5). Considering $\psi
(z/\vert \psi' (0)\vert)$ and $\chi (z/\vert \chi'(0)\vert)$ we find
$r, s
> 0$ such that the following holds:

\begin{itemize}
\item[$\alpha$)] There are functions from ${\mathcal O}
{\mathbf L}_0^*$, denoted again by $\psi$ and $\chi$, which are
defined on $B_{{\mathbf L}}(r)$ and fulfil $\vert a (\psi) \vert =
\vert a (\chi) \vert = 1$, such that $\partial D \cap B_{{\mathbf
L}}(s) = \Gamma \cup \Gamma'$ with $\Gamma := \psi (] 0, \varepsilon
])$ and $\Gamma' := \chi (] 0, \varepsilon[)$ for some $0 <
\varepsilon \leq r$. Moreover, we can assume that $\psi^{-1}$ and
$\chi^{-1}$ are defined on $B_{{\mathbf L}}(r)$ and that $\psi,
\chi, \psi^{-1}, \chi^{-1}$ are injective on $B_{{\mathbf L}}(r)$.
\item[$\beta$)] There are Puiseux series $h_0, h_1$ convergent on
$B_{{\mathbf L}}(s)$ such that $h = h_0$ on $\Gamma$ and $h = h_1$
on $\Gamma'$.
\end{itemize}

\noi We may assume that $s \leq r$ and that $\arg (a(\psi)) < \arg
(a(\chi))$. Recursively we find for $k \geq 1$ constants $r_k \geq
s_k > 0$, a domain $D_k \subset {\mathbf L}$, a function $f_k \in
{\mathcal O} (D_k)$ and a function $\varphi_k \in {\mathcal O}
{\mathbf L}_0^*$ with $\vert a (\varphi_k) \vert = 1$, $r_k = r
(\varphi_k) = r (\varphi_k^{-1})$ and a Puiseux series $h_k$
convergent on $B_{{\mathbf L}}(s_k)$ with the following properties:
\begin{itemize}
\item[a)] $D_k \supset D_{k-1}$ and $f_k$ extends $f_{k-1}$ holomorphically for $k
\geq 2$.
\item[b)] $\partial D_k \cap B_{{\mathbf L}} (s_k)$ consists of two boundary
curves $\Gamma_k$ and $\Gamma'_k$ where $\Gamma_k := \psi (]0,
\mu_k[) \subset \Gamma$ and $\Gamma'_k := \varphi_k (]0,
\varepsilon_k[)$ for some $0 < \mu_k$, $\varepsilon_k \leq r_k$,
\item[c)] $\Re\, f_k = h_k$ on $\Gamma_k$
\end{itemize}

\noi as follows:\\
$k = 1$: $r_1 := r$, $s_1 := s$, $\varphi_1 := \chi$, $D_1 := D$,
$f_1 := f$ from Theorem 2.9, and $h_1$ from $\beta$) above.\\
$k \to k+1$: We set $r_{k+1} := \frac{r_k}{100}$, $s_{k+1} :=
\frac{s_k}{100}$, $D_{k+1} := (\overline{D_k \cup D'_k})^0$ with
$$D'_k := \varphi_k (\tau (\varphi_k^{-1} (D_k \cap B\!\left(
\frac{s_k}{4} \right))),$$

\medskip\noi
$(\dagger) \qquad f_{k+1} \colon D_{k+1} \to {\mathbb C},$

\noindent
$\qquad\qquad\qquad z \mapsto \left\{ \begin{array}{lll} f_k (z) & & z \in D_k, \\
         & \mbox{if} &       \\
- \overline{((f_k - h_k) \circ \varphi_k \circ \tau \circ
\varphi_k^{-1} (z))}
+ h_k (z) & & z \in D'_k, \\
\end{array} \right. $

\medskip
\noi $\varphi_{k+1} := \varphi_k \circ \tau \circ \varphi_k^{-1}
\circ \psi \circ \tau$ and $h_{k+1} := - \overline{(h_0 - h_k) \circ
\varphi_k \circ \tau \circ \varphi_k^{-1}} + h_k$.

\medskip\noi
Note that these data are well defined:\\
By Proposition 1.9 a) $D'_k$ is well defined and $f_{k+1}$ exists on
$D_{k+1}$. Note that the angle of $D_{k+1}$ at 0 is twice the angle
of $D_k$ at 0. The function $f_{k+1}$ originates from a reflection
process at $\Gamma_k$ and extends $f_k$ holomorphically (compare
with Example 3.4). By Proposition 1.9 b) and Remark 3.5 we can
choose $r(\varphi_{k+1}) = r (\varphi_{k+1}^{-1}) =
\frac{r_k}{100}$. By Proposition 1.9 a), Proposition 1.10 and Remark
3.5 $h_{k+1}$ is a Puiseux series convergent on $B_{{\mathbf
L}}\!\left( \frac{s_k}{4} \right)$. By construction and by
Proposition 1.9 a) we see
that a), b) and c) hold.\\
We extend $f$ holomorphically to $\bigcup\limits_k D_k$ by setting
$f \vert_{D_k} := f_k$. We see that $s_k = s/100^{k-1}$ and $\arg (a
(\varphi_k)) - \arg (a(\psi)) = 2^{k-1} (\sphericalangle D)$. By the
definition of $\Gamma_k$ and Proposition 1.9 a) we get that
$$\Gamma_k \subset \left\{ (r, \varphi) \in {\mathbf L} \; \vert \;
\vert \varphi - \arg a (\varphi_k) \vert \leq \frac{\pi}{2}
\right\}.$$

\noi Hence we see by b) that $f$ is holomorphic on
$$\left\{
(r, \varphi) \in {\mathbf L} \; \vert \; 0 < r < s/100^{k-1} \;
\mbox{and} \; 0 < \varphi - \alpha < 2^{k-1} (\sphericalangle D) -
\frac{\pi}{2} \right\}$$

\noi for all sufficiently large $k \geq 1$ where $\alpha := \arg (a
(\psi))$. For $\varphi > \alpha$ let $k (\varphi) \in {\mathbb N}$
be such that $2^{k-1} (\sphericalangle D) - \frac{\pi}{2} \leq
\varphi - \alpha \leq 2^k (\sphericalangle D) - \frac{\pi}{2}$. Then
there is some $C
> 0$ such that $k (\varphi) \leq C \log (\varphi - \alpha)$ for
all sufficiently large $\varphi > 0$. Hence we find some $K
> 1$ such that $f$ is holomorphic on
$$\left\{ (r, \varphi) \in {\mathbf L} \; \vert \; \varphi > \alpha
\quad \mbox{and} \quad 0 < r < K^{- \log_+ (\varphi - \alpha)}
\right\}$$

\noi where $\log_+ x := \max \{ 1, \log x \}$ for $x > 0$.

\medskip
\noi Repeating this process in the negative direction we see that
$f$ is holomorphic on some quadratic domain $U$ since $\log
(\varphi- \alpha) \leq \sqrt{\varphi}$ for all sufficiently large
$\varphi$.

\medskip\noi
{\bf Step 2:} We show that $f \in {\mathcal Q}^{\log} (U)$. Let $g
:= \sum\limits_{\alpha \geq 0} a_{\alpha} P_{\alpha} (\log z)
z^{\alpha} \in {\mathbb C} [[z^*]]_{\log}$ from Theorem 2.9. Given
$R > 0$ we show that there is a quadratic domain $U_R$ such that
$$f(z) - \sum\limits_{\alpha \leq R} a_{\alpha} P_{\alpha} (\log z)
z^{\alpha} = o (\vert z \vert^R) \quad \mbox{as} \quad \vert z \vert
\to 0 \quad \mbox{on} \quad U_R.$$

\noi Hence $Tf = g$ and we see with Theorem 2.9 that $f \in
{\mathcal Q} (U)$ if $\sphericalangle D / \pi \in {\mathbb R}
\setminus {\mathbb Q}$.

\bigskip\noi
Proof of Step 2:\\
We work in the setting of Step 1. We can assume that $h_0$ and $h_1$
have a common denominator $d \in {\mathbb N}$ and that there is a $c
> 0$ such that $\vert h_i (z) \vert \leq c$ for $\vert z \vert < s$,
$i = 0, 1$. By Proposition 1.9 a), Proposition 1.10 and Remark 3.5
we see that $h_k$ is a Puiseux series on $B_{{\mathbf L}}\!\left(
\frac{s_{k-1}}{4} \right)$ with denominator $d$ for all $k \geq 2$.
Moreover, by induction we see that $\vert h_k (z) \vert \leq 3^{k-1}
c$ for $\vert z \vert <
\frac{s_{k-1}}{4}$ and all $k \geq 2$.\\
We fix $R > 0$. Let $\gamma := \sum\limits_{\alpha \leq R}
a_{\alpha} P_{\alpha} (\log z) z^{\alpha} \in {\mathbb C}
[[z^*]]_{\log}^{\fin}$. We choose $R'
> R$ with the following properties:

\medskip\noi
(i) $R' < \min \{ \alpha \in \supp (g) \; \vert \; \alpha > R \}$,\\
(ii) $R' < \frac{[Rd]+1}{d}$,\\
(iii) $R' < [R - \alpha] + \alpha + 1$ for all $\alpha \in \supp
(\gamma)$.

\medskip
\noi Here $[x]$ denotes the largest integer $n \leq x$. Let $K$ be
the number of elements of $\supp (\gamma)$, $L := \sup \{ \vert a
\vert : a$ is a coefficient of $a_{\alpha} P_{\alpha}$, $\alpha \in
\supp (\gamma) \}$ and $M := \max\limits_{\alpha \in \supp (\gamma)}
\deg
P_{\alpha}$.\\
We may choose $r$ in Step 1 such that the following holds:

\medskip\noi (iv) $\vert z^{\alpha + N} (\log z)^M \vert \leq \vert z
\vert^{R'}$ on $B_{{\mathbf L}}(r)$ for all $N > [R - \alpha]$ and
all $\alpha \in \supp (\gamma)$.

\medskip\noi By Theorem 2.9 and condition (i) we find $C > 1$ such that $\vert f(z) -
\gamma (z) \vert \leq C \vert z \vert^{R'}$ on $\overline{D}$.

\bigskip \noi {\bf Claim 1:} Given $k \geq 1$ there is some $C_k
> 0$ such that
$$\vert f_k (z) - \gamma (z) \vert \leq C_k \vert z
\vert^{R'}$$

\noindent on $\overline{D_k} \cap B_{{\mathbf L}} (s_k)$.

\medskip\noindent
From the proof of Claim 1 we will obtain below estimates for the
$C_k$'s.

\bigskip\noi
Proof of Claim 1 by induction on $k$:\\
The case $k = 1$ is settled by the above.

\noi $k - 1 \to k$: By construction we have that $f_{k}
\vert_{D_{k-1} \cap B_{{\mathbf L}} (s_k)} = f_{k-1} \vert_{D_{k-1}
\cap B_{{\mathbf L}} (s_k)}$. For the inductive step we have to
consider $f_{k} \vert_{D'_{k} \cap B_{{\mathbf L}}(s_{k})}$. From
($\dagger$) in Step 1 we see that
$$f_{k} \vert_{D'_{k} \cap B_{{\mathbf L}} (s_{k})} =
v_{k, 1} + v_{k, 2} + v_{k, 3}$$

\noi with $v_{k, 1} := - \overline{f_{k-1} \circ \varphi_{k-1} \circ
\tau \circ \varphi_{k-1}^{-1}}, v_{k, 2} := \overline{h_{k-1} \circ
\varphi_{k-1} \circ \tau \circ \varphi_{k-1}^{-1}}$ and $v_{k, 3} :=
h_{k-1}$. To prove Claim 1 we show the following

\bigskip
\noi {\bf Claim 2:} Let $1 \leq i \leq 3$. Assuming that Claim 1
holds for $k-1$ we show that there is some $\gamma_{k, i} \in
{\mathbb C} [[z^*]]_{\log}^{\fin}$ with $\supp (\gamma_{k, i})
\subseteq [0, R]$ and some $C_{k, i} > 0$ such that
$$\vert v_{k, i} - \gamma_{k, i} \vert \leq C_{k, i} \vert z \vert^{R'}$$

\noi on $B_{{\mathbf L}} (s_k)$.

\medskip
\noi Claim 1 follows from Claim 2: since $f_{k} = f_{k-1}$ on
$\Gamma_{k-1}$ we see by applying the inductive hypothesis and by
Claim 2 that $\gamma = \gamma_{k,1} + \gamma_{k,2} + \gamma_{k, 3}$
and obtain that Claim 1 holds for $k$ with $C_{k} := C_{k, 1} +
C_{k, 2} + C_{k, 3}$. \hfill (1)

\bigskip\noi Proof of Claim 2:\\
To prove Claim 2 we distinguish the three cases $i = 1, 2, 3$.\\
$i = 1$: Let $\delta_{k-1} := f_{k-1} - \gamma$. Then $\vert
\delta_{k-1} (z) \vert \leq C_{k-1} \vert z \vert^{R'}$ for $z \in
D_{k-1} \cap B_{{\mathbf L}} (s_{k-1})$ by the assumption of Claim
2. Let $\eta_{k} := - \overline{\delta_{k-1} \circ \varphi_{k-1}
\circ \tau \circ \varphi_{k-1}^{-1}}$. Then $v_{k, 1} = -
\overline{\gamma \circ \varphi_{k-1} \circ \tau \circ
\varphi_{k-1}^{-1}} + \eta_{k}$ and $\vert \eta_k (z) \vert \leq
C_{k-1} 4^{R'} \vert z \vert^{R'}$ on $D'_{k} \cap B_{{\mathbf L}}
(s_{k})$ by Proposition 1.9 a). Let $w_{k, 1} := \gamma \circ
\varphi_{k-1}$, $w_{k, 2} := - \overline{w_{k,1} \circ \tau}$ and
$w_{k, 3} : = w_{k, 2} \circ \varphi_{k-1}^{-1}$. To prove Claim 2
in the case $i = 1$ we show the following

\bigskip\noi
{\bf Claim 3:} Let $1 \leq j \leq 3$. We show that there is some
$\lambda_{k, j} \in {\mathbb C} [[z^*]]_{\log}^{\fin}$ with $\supp
(\lambda_{k, j}) \subseteq [0, R]$ and some $E_{k, j} > 0$ such that
$$\vert
w_{k, j} - \lambda_{k, j} \vert \leq E_{k, j} \vert z \vert^{R'}
$$

\noi on $B_{{\mathbf L}}(s_{k})$.

\medskip\noi
The case $i = 1$ from Claim 2 follows from Claim 3: assuming Claim
3, we can choose $\gamma_{k, 1} := \lambda_{k, 3}$ and $C_{k, 1} : =
E_{k, 3} + C_{k-1} 4^{R'}$. \hfill (2)

\bigskip\noi Proof of Claim 3:\\
To prove Claim 3 we distinguish the three cases
$j = 1, 2, 3$.\\
$j = 1$: Let $\alpha \in \supp (\gamma)$ and $m \leq M$. By Lemma
1.20 we find power series $p_1, \dots, p_m$ convergent on
$B_{{\mathbf L}} (r_{k-1})$ such that
$$(z^{\alpha} (\log z)^m) \circ
\varphi_{k-1} = z^{\alpha} \left( \sum\limits_{\ell = 0}^m p_{\ell}
(z) (\log z)^{\ell} \right)$$

\noindent with $\vert p_{\ell} (z) \vert \leq 2^{m+\alpha} (\vert
\arg a (\varphi_{k-1}) \vert + 3)^m$. Hence
$$v_{k, 1} = \gamma \circ \varphi_{k-1} = \sum\limits_{\alpha \in \supp (\gamma)}
z^{\alpha} \left( \sum\limits_{\ell = 0}^M q_{\alpha, \ell} (z)
(\log z)^{\ell} \right)$$

\noi with power series $q_{\alpha, \ell}$ convergent on $B_{{\mathbf
L}} (r_{k-1})$ and $\vert q_{\alpha, \ell} \vert \leq d_{k}$ on
$B_{{\mathbf L}}(r_{k-1})$ where $d_{k} := KL 2^{M+R} (\vert \arg a
(\varphi_{k-1}) \vert + 3)^M$. Let $T_{\alpha, \ell}$ be the
truncated power series expansion of $q_{\alpha, \ell}$ up to order
$[R - \alpha]$. By Proposition 1.5 b) we get $q_{\alpha, \ell} (z) -
T_{\alpha, \ell} (z) \vert \leq \frac{2 d_k}{r_{k-1}^{R+1}} \vert z
\vert^{[R- \alpha ] + 1}$ on $B_{{\mathbf L}}\!\left(
\frac{r_{k-1}}{2} \right)$. By condition (iii) and (iv) we see that
there is $\lambda_{k, 1} = \sum\limits_{\alpha \leq R} b_{\alpha}
Q_{\alpha} (\log z) z^{\alpha} \in {\mathbb C}
[[z^*]]_{\log}^{\fin}$ such that on $B_{{\mathbf L}} \left(
\frac{r_{k-1}}{2} \right)$
$$\vert w_{k, 1} (z)  - \lambda_{k, 1} (z) \vert \leq E_{k, 1} \vert z \vert^{R'}$$

\noi where $E_{k, 1} := \frac{2d_{k}}{r_{k-1}^{R+1}}$. \hfill (3)

\noindent Moreover, let $K'$ be the number of elements of $\supp
(\lambda_{k, 1})$, $L' := \sup \{ \vert a \vert : a$ is a
coefficient of $b_{\alpha} Q_{\alpha}$, $\alpha \in \supp
(\lambda_{k, 1}) \}$ and $M' := \max\limits_{\alpha \in \supp
(\lambda_{k, 1})} \deg Q_{\alpha}$. Then $K' \leq K + R$, $L' \leq
\frac{d_{k}}{r_{k-1}^{R}} L$ by Cauchy estimates and $M' \leq M$.

\medskip\noi
$j = 2$: We set $\lambda_{k, 2} := - \lambda_{k, 1}$ and $E_{k, 2}
:= E_{k, 1}$. \hfill (4)

\noindent By Remark 3.5 we see that on $B_{{\mathbf
L}}\!(\frac{r_{k-1}}{2})$
$$\vert w_{k, 2} (z) - \lambda_{k, 2} (z) \vert \leq E_{k, 2} \vert z
\vert^{R'}. \leqno(\ast)$$

\medskip\noi
$j = 3$: As in the case $j = 1$ we see that there is $\lambda_{k, 3}
\in {\mathbb C} [[z^*]]_{\log}^{\fin}$ with $\supp (\lambda_{k, 3})
\subseteq [0, R]$ such that
$$\vert \lambda_{k, 2} \circ \varphi_k^{-1} (z) - \lambda_{k, 3} (z) \vert \leq
\frac{2 d'_{k}}{r_{k-1}^{R}} \vert z \vert^{R'} \leqno(\ast \ast)$$

\noi on $B_{{\mathbf L}}\!\left( \frac{r_{k-1}}{2} \right)$ with
$d'_{k} := K' L' 2^{M+R} [\vert \arg a (\varphi_{k-1}) \vert + 3]^M$
(note that $\arg a (\varphi_{k-1}^{-1}) = - \arg a
(\varphi_{k-1})$). Hence we obtain by applying Proposition 1.9 a) to
($\ast$) and by ($\ast \ast$) that on $B_{{\mathbf L}}\!\left(
\frac{r_{k-1}}{4} \right)$
$$\vert w_{k, 3} (z) -
\lambda_{k, 3} (z) \vert \leq E_{k, 3} \vert z \vert^{R'}
$$

\noi where
$$E_{k, 3} := \frac{2 d'_{k}}{r_{k-1}^{R'}} + 2^{R'} E_{k, 2}.\eqno(5)$$

\noi So Claim 3 and therefore the case $i = 1$ of Claim 2 is proven
since $s_k < \frac{r_{k-1}}{4}$. We continue with the case $i = 2$
of Claim 2.

\medskip\noi
$i = 2$: We see with Proposition 1.9 a), Proposition 1.10 and Remark
3.5 that $v_{k, 2}$ is a Puiseux series on $B_{{\mathbf L}}\!\left(
\frac{s_{k-1}}{4} \right)$ with denominator $d$. Moreover, $\vert
v_{k, 2} \vert \leq 3^{k-1} c$ on $B_{{\mathbf L}}\!\left(
\frac{s_{k-1}}{4} \right)$ (compare with the beginning of the proof
of Step 2). Hence by Proposition 1.5 b) and condition (ii) we find
$\gamma_{k, 2}$ as described such that $\vert v_{k, 2} (z) -
\gamma_{k, 2} (z) \vert \leq C_{k, 2} \vert z \vert^{R'}$ on
$B_{{\mathbf L}} (s_{k})$ with
$$C_{k, 2} := 3^{k-1} c d'
\left(\frac{4}{s_{k-1}} \right)^{\frac{R'+1}{d}} \eqno(6)$$

\noindent and $d' := \frac{1}{1- (\frac{1}{25})^{\frac{1}{d}}}$.

\medskip\noi
$i = 3$: As in the case $i = 2$ we can find with Proposition 1.5 b)
$\gamma_{k, 3}$ and $C_{k, 3}$ as described. Moreover, we can choose
$C_{k, 3} := C_{k, 2}$. \hfill (7)

\medskip\noi
Hence Claim 2 and therefore also Claim 1 is proven. We revisit now
the construction of the constants $C_k$ of Claim 1. By (1) we have
$C_k = C_{k, 1} + C_{k, 2} + C_{k, 3}$ with $C_{k, i}$ from Claim 2.
By (2) we have $C_{k, 1} = E_{k, 3} + C_{k-1} 4^{R'}$ with $E_{k,
3}$ from Claim 3. The constant $E_{k, 3}$ is computed via (3), (4)
and (5). Since $r_{k-1} = \frac{r}{100^{k-2}}$ and $\vert \arg a
(\varphi_{k-1}) - \arg a (\psi) \vert = 2^{k-1} (\sphericalangle D)$
we see that there is some $B > 1$ such that $E_{k, 3} \leq B^k$ for
all $k \geq 2$. Since $s_{k-1} = \frac{s}{100^{k-2}}$ we see by (6)
that after enlarging $B$ if necessary  $C_{k, 2} \leq B^k$ for all
$k \geq 2$. Since $C_{k, 2} = C_{k, 3}$ by (7) we finally find by
(1) and (2) some $A > 1$ such that $C_k \leq A^k$ for all $k$. Hence
Claim 1 gives that $\vert f (z) - \gamma (z) \vert \leq A^{k} \vert
z \vert^{R'}$ on $D_k \cap B_{{\mathbf L}} (s_k)$ and all $k \in
{\mathbb N}$. Moreover, we may assume that $A^{k} \leq s_k$ for all
$k \in {\mathbb N}$. We choose $R < S < R'$ with $R' - S \leq 1$. We
set $t_k := A^{- \frac{k}{R'-S}}$ and obtain on $D_k \cap
B_{{\mathbf L}}(t_k)$
$$\vert f(z) - \gamma (z) \vert \leq A^{k} \vert z \vert^{R'-S} \vert
z \vert^S \leq \vert z \vert^S.$$

\noi Using a similar argument as at the end of the proof of Step 1
we find some $K > 1$ such that $\vert f(z) - g_0 (z) \vert \leq
\vert z \vert^S$ on the set
$$\{ (r, \varphi) \in {\mathbf L} \; \vert \; \varphi > \alpha
\quad \mbox{and} \quad 0 < r < K^{- (\log_+ (\varphi- \alpha))} \}$$

\noi where $\alpha := \arg (a (\psi))$ and $\log_+ x := \max \{ 1,
\log x \}$ for $x > 0$. Repeating these arguments in the negative
direction we see that $f(z) - \gamma (z) = o ( \vert z \vert^R)$ on
some admissible domain $U_R \subset U$ since $S > R$ and $\log
(\varphi -\alpha) \leq \sqrt{\varphi}$ eventually. \hfill $\square$

\bigskip\noi We need the following final ingredients for the proof of
Theorem B.

\bigskip\noi
{\it Definition 3.6} (compare with [27, Definition 3.4 \& Definition
4.3])

\noi Let $\lambda \in \overline{{\mathbb H}} \setminus \{ 0 \}$
where ${\mathbb H}$ denotes the upper half plane. We have $B( \vert
\lambda \vert, \vert \lambda \vert) \subset {\mathbf L}$ via the
identification of ${\mathbb C} \setminus {\mathbb R}_{\leq 0}$ with
${\mathbb R}_{> 0} \times ] - \pi, \pi [ \subset {\mathbf L}$. Let
$\lambda = \vert \lambda \vert \, e^{ia}$ with $0 \leq a \leq \pi$.
We identify $B(\lambda, \vert \lambda \vert)$ with $\{ (r, \varphi)
\in {\mathbf L} \; \vert \; (r, \varphi - a) \in B (\vert \lambda
\vert, \vert \lambda \vert) \}$. We set ${\mathbf t}_{\lambda}
\colon B (0, \vert \lambda \vert) \to B (\lambda, \vert \lambda
\vert)$, $z \mapsto \lambda + z$, and for $\rho
> 0$ we define ${\mathbf r}^{\rho, \lambda} \colon {\mathbf L}
\times B (0, \vert \lambda \vert) \to {\mathbf L}^2$, $(z_1, z_2)
\mapsto (z_1, w_2)$, where $w_2 : = {\mathbf m} ({\mathbf p}^{\rho}
(z_1), {\mathbf t}_{\lambda} (z_2))$.

\bigskip\noi
{\it Remark 3.7}

\noi Let $U \subset {\mathbf L}^2$ be a 2-quadratic domain (compare
with [27, Definition 2.4]) and let $f \in {\mathcal Q}_2^2 (U)$
(compare with [27, Definition 5.1]). Let $\lambda \in
\overline{{\mathbb H}} \setminus \{ 0 \}$. As in [27,
Proposition~4.4 \& Proposition~5.15] we find some 1-quadratic domain
$V \subset {\mathbf L} \times B_{{\mathbf L}} (\lambda)$ such that
${\mathbf r}^{1, \lambda} (V) \subset U$ and the function ${\mathbf
r}^{1, \lambda} f : = f \circ {\mathbf r}^{1, \lambda} \in {\mathcal
Q}_1^2 (V)$.

\bigskip\noi
PROPOSITION 3.8

\noi {\it Let $U \subset {\mathbf L}$ be a quadratic domain and let
$f \in {\mathcal Q} (U)$. Let $V \subset U$ be a domain such that
$e$ is injective on $V$. Then there is some $r
> 0$ such that $f \circ (e \vert_V)^{-1} \vert_{B(0, r)}$ is
definable in} ${\mathbb R}_{{\mathcal Q}}$.

\bigskip\noi {\it Proof}

\noi By considering finite coverings we can assume that $V := \{ (r,
\varphi) \in U \; \vert \; \varphi_0 < \varphi < \varphi_0 + \pi \}$
for some $\varphi_0 \in {\mathbb R}$. Let $b := (1, \varphi_0) \in
{\mathbf L}$ and $\psi := {\mathbf m} (b, {\mathbf p}^1 (z)) \in
{\mathcal O} {\mathbf L}_0^*$. Then $f \circ \psi \in {\mathcal Q}
(W)$ for some quadratic domain $W \subset {\mathbf L}$ by
Proposition 1.19 c). Let $V' := \{ (r, \varphi) \in W \; \vert \; 0
< \varphi < \pi \}$. Then $f \circ (e \vert_V)^{-1} \vert_{B(0, r)}
= f \circ \psi \circ (e \vert_{V'})^{-1} \vert_{B(0, r)}$ for all
sufficiently small $r > 0$; we can therefore assume that $\varphi_0
= 0$. Hence via the identification of ${\mathbb C} \setminus
{\mathbb R}_{\leq 0}$ with $\{ (r, \varphi) \in {\mathbf L} \; \vert
\; - \pi < \varphi < \pi \}$, it suffices to show that $f
\vert_{{\mathbb H} \cap B (0, r)}$ is definable in ${\mathbb
R}_{{\mathcal Q}}$ for some $r
> 0$. We define $\Phi \colon U \times U \to {\mathbb C}$, $(z_1,
z_2) \mapsto f(z_2)$. Let $a \in [0, \pi]$. We consider $g_a :=
{\mathbf r}^{1, \lambda_a} \Phi$ with $\lambda_a := e^{ia}$. By
Remark 3.7 we have that $g_a \in {\mathcal Q}_1^2$. We set $G_a :=
g_a (z_1, h_a (z_2))$ with $h_a (z) := e^{i (z+a)} - e^{ia}$. Then
$G_a \in {\mathcal Q}_1^2$ by [27, Proposition 5.10]. Hence there is
some $r_a > 0$ and some quadratic domain $U_a$ such that $G_a \in
{\mathcal Q}_1^2 (U_a \times B(0, r_a))$. We can assume that $U_a =
\{ (r, \varphi) \in {\mathbf L} \; \vert \; 0 < r < c_a \exp (- C_a
\sqrt{\vert \varphi \vert}) \}$ with some positive constants $c_a$,
$C_a$. We define $G^*_a \colon U_a \times B(0, r_a) \to {\mathbb
C}$, $(z_1, z_2) \mapsto \overline{G_a (\tau (z_1),
\overline{z_2})}$. Note that $G^*_a \in {\mathcal Q}_1^2$ (compare
with the proof of [27, Proposition 7.3]. We set ${\mathcal R} G_a :
= \frac{1}{2} (G_a + G^*_a)$ and ${\mathcal J} G_a : = \frac{1}{2i}
(G_a - G^*_a)$. Then ${\mathcal R} G_a$, ${\mathcal J} G_a \in
{\mathcal Q}_{1, 1; \varepsilon_a}$ for some $\varepsilon_a > 0$
(compare with [27, section 7]. Hence ${\mathcal R} G_a$ and
${\mathcal J} G_a$ are defined on $I_a : = [ 0, \varepsilon_a]
\times [- \varepsilon_a, \varepsilon_a]$, and ${\mathcal R} G_a
\vert_{I_a}$ and ${\mathcal J} G_a \vert_{I_a}$ are definable in
${\mathbb R}_{{\mathcal Q}}$.\\
For $(r, \varphi) \in I_a$ we have ${\mathcal R} G_a (r, \varphi) =
\Re f ( r e^{i (\varphi + a)}) = \Re f (r \cos (\varphi + a), r \sin
(\varphi + a))$ and ${\mathcal J} G_a (r, \varphi) = \Im f (r e^{i
(\varphi + a)}) = \Im f (r \cos (\varphi + a), r \sin (\varphi +
a))$. Since the polar coordinates are definable in ${\mathbb
R}_{{\mathcal Q}}$ we find by a compactness argument (note that $a
\in [0, \pi ]$) some $r > 0$, such that $f \vert_{{\mathbb H} \cap B
(0, r)}$ is definable in ${\mathbb R}_{{\mathcal Q}}$. \hfill
$\square$

\bigskip\noi We obtain Theorem B:

\bigskip\noi THEOREM 3.9

\noi {\it Let $\Omega \subset {\mathbb R}^2$ be a semianalytic and
bounded domain without isolated boundary points. Suppose that}
$\sphericalangle (\Omega, x) \subset \pi ({\mathbb R} \setminus
{\mathbb Q})$ {\it for all} $x \in \Sing (\partial \Omega)$. {\it
Let} $h \colon
\partial \Omega \to {\mathbb R}$ {\it be a continuous and semianalytic
function on the boundary and let $u$ be the Dirichlet solution for
$h$. Then $u$ is definable in the o-minimal structure} ${\mathbb
R}_{{\mathcal Q}, \exp}$.

\bigskip\noi {\it Proof}

\noi Let $x \in \overline{\Omega}$. We need to show that there is a
neighbourhood $V_x$ of $x$ such that $u \vert_{V_x \cap \Omega}$ is
definable in ${\mathbb R}_{{\mathcal Q}, \exp}$.

\medskip\noi
Case 1: $x \in \Omega$. Let $r > 0$ with $B(x, r) \subset \Omega$.
Then $u \vert_{B(x, r)}$ is harmonic and therefore real analytic
(see [2, Theorem 1.8.5]). Hence $u \vert_{B \left(x, \frac{r}{2}
\right)}$ is definable in ${\mathbb R}_{\an}$ which is a reduct of
${\mathbb R}_{{\mathcal Q}, \exp}$.

\medskip\noi
Case 2: $x \in \partial \Omega$

\medskip\noi
Case 2.1: $x \not\in \Sing (\partial \Omega)$. Then $\partial
\Omega$ is a real analytic manifold at $x$. Let $r > 0$ such that
$\Omega' := \Omega \cap B (x, r)$ is simply connected. Let
$\tilde{h} := u \vert_{\partial \Omega'}$. By the Riemann Mapping
Theorem there is a biholomorphic map $\varphi \colon \Omega' \to
B(0, 1)$. By Caratheodory's theory of prime ends (see Pommerenke
[33, Chapter 2, p.18]) and by the curve selection lemma (see [10,
p.94]) $\varphi$ has a continuous extension to $x$. Then by the
Schwarz reflection principle $\varphi$ has a holomorphic extension
to a neighbourhood of $x$. Let $\widehat{h} := \tilde{h} \circ
\varphi$ and $y := \varphi (x)$. Then there is some $s > 0$ such
that $\widehat{h} \vert_{B (y, s)}$ is semianalytic. Let
$\widehat{u}$ be the Dirichlet solution for $\widehat{h}$ on $B(0,
1)$. Then $\widehat{u}$ is given by the Poisson integral (see [20,
Theorem 1.8]):
$$\widehat{u} (\xi) = \frac{1}{2 \pi} \int\limits_{\partial B (0,1)}
\frac{1 - \vert \xi \vert^2}{\vert \eta - \xi \vert^2} \widehat{h}
(\eta) d \sigma (\eta).$$

\noi We have therefore
$$\begin{array}{rcl}
\widehat{u} (\xi) & = & \frac{1}{2 \pi} \int\limits_{\partial B
(0,1) \setminus B (y, s)} \frac{1 - \vert \xi \vert^2}{\vert \eta -
\xi \vert^2} \widehat{h} (\eta) d \sigma (\eta) \\
& + & \frac{1}{2 \pi} \int\limits_{\partial B (0,1) \cap B (y, s)}
\frac{1 - \vert \xi \vert^2}{\vert \eta -
\xi \vert^2} \widehat{h} (\eta) d \sigma (\eta). \\
\end{array}$$

\noi Applying the Laplace operator to the first integral we see by
switching differentation and integration that the first summand is
harmonic in $B(y, s)$ and hence real analytic in $B(y, s)$; the
second one is definable in the o-minimal structure ${\mathbb
R}_{\an, \exp}$ by [7, Th\'eor\`eme 1$'$]. Hence $\widehat{u}
\vert_{B \left( y, \frac{s}{2} \right)}$ is definable in ${\mathbb
R}_{\an, \exp}$. Let $0 < t < r$ be such that $\varphi$ has a
holomorphic extension to $B(x, 2t)$ and $\varphi (B(x, t)) \subset B
\left( y, \frac{s}{2} \right)$. Then $u \vert_{B(x, t)} =
\widehat{u} \circ \varphi \vert_{B(x, t)}$ is definable in ${\mathbb
R}_{\an, \exp}$ which is a reduct of ${\mathbb R}_{{\mathcal Q},
\exp}$.

\medskip\noi
Case 2.2: $x \in \Sing (\partial \Omega)$. We may assume that $x =
0$. There is a semianalytic neighbourhood $V$ of 0 such that $V \cap
\Omega$ is the disjoint union of finitely many corner components of
$\Omega$ at 0. Let $C$ be such a corner component. Then
$\sphericalangle_0 C \in \pi ({\mathbb R} \setminus {\mathbb Q})$ by
assumption. By Theorem 3.2 there is a quadratic domain $U \subset
{\mathbf L}$ and an $f \in {\mathcal Q} (U)$ such that $\Re f$
extends $u \vert_C$. With Proposition 3.8 we get that $u \vert_C$ is
definable in ${\mathbb R}_{{\mathcal Q}}$ and hence in ${\mathbb
R}_{{\mathcal Q}, \exp}$. This shows the claim. \hfill $\square$

\bigskip\noi
As an application we obtain the definability of the Green function.

\bigskip\noi
COROLLARY 3.12

\noi {\it Let} $\Omega \subset {\mathbb R}^2$ {\it be a semianalytic
and bounded domain without isolated boundary points. Suppose that}
$\sphericalangle (\Omega, x) \subset \pi ({\mathbb R} \setminus
{\mathbb Q})$ {\it for all} $x \in \Sing (\partial \Omega)$. {\it
Let $y \in P$. Then the Green function of $\Omega$ with pole $y$ is
definable in} ${\mathbb R}_{{\mathcal Q}, \exp}$.

\bigskip\noi
{\it Proof}

\noi Let $y \in P$ and $K_y \colon {\mathbb R}^2 \setminus \{ y \}
\to {\mathbb R}$, $x \mapsto \log \frac{1}{\vert x - y \vert}$, be
the Poisson kernel with pole $y$. Then the Green function of
$\Omega$ with pole $y$, denoted by $G_y^{\Omega}$, is given by
$G_y^{\Omega} = K_y - u$, where $u$ is the Dirichlet solution for
$K_y \vert_{\partial \Omega}$ (see [2, Chapter 4.1]). Since $K_y
\vert_{\partial \Omega}$ is semianalytic we have by Theorem 3.9 that
$u$ is definable in ${\mathbb R}_{{\mathcal Q}, \exp}$. Since $K_y$
is definable in ${\mathbb R}_{\exp}$ we obtain the claim. \hfill
$\square$

\bigskip\noi
In the case of semilinear domains we can overcome the restriction on
the angles:

\bigskip\noi
COROLLARY 3.13

\noi {\it Let $\Omega \subset {\mathbb R}^2$ be a semilinear domain
without isolated boundary points. Let $h \colon \partial \Omega \to
{\mathbb R}$ be a continuous and semianalytic function on the
boundary and let $u$ be the Dirichlet solution for $h$. Then $u$ is
definable in the o-minimal structure} ${\mathbb R}_{{\mathcal Q},
\exp}$.

\bigskip
\noi {\it Proof}

\noi Let $x \in \Sing (\partial \Omega)$ and let $C$ be a semilinear
representative of a component of the germ of $\Omega$ at $x$ such
that the germ of $C$ at $x$ is connected. If $\sphericalangle_x C /
\pi \in {\mathbb Q}$ we see with [24, Corollary 4] that $u \vert_C$
is definable in the o-minimal structure ${\mathbb R}_{\an, \exp}$.
With the proof of Theorem 3.9 we obtain the claim. \hfill $\square$

\vspace{0.5cm}\noi {\small

}

\vspace{0.5cm}\noindent {\small University of Regensburg, Department
of Mathematics, Universit\"atsstr. 31,\\ D-93040 Regensburg,
Germany; Tobias.Kaiser@mathematik.uni-regensburg.de}

\end{document}